\documentclass[letterpaper, 11pt,  reqno]{amsart}
\usepackage[margin=1.2in,marginparwidth=1.5cm,
 marginparsep=0.5cm]{geometry}
\usepackage{amsmath,amssymb,amscd,amsthm,amsxtra,
  esint,xcolor,mathtools,mathrsfs}

\usepackage[implicit=true,colorlinks=true,
  linkcolor=blue,citecolor=blue,urlcolor=blue]{hyperref}
\definecolor{labelkey}{rgb}{0,0,1}
\usepackage{graphicx}

\usepackage{bm,cases}
\usepackage{enumitem}
\usepackage{booktabs,tabularx}
\setlist[enumerate]{itemsep=.5em,label=\normalfont(\roman*)}
\usepackage[norefs,nocites]{refcheck}

\definecolor{gr}{rgb}   {0.,   0.69,   0.23 }
\definecolor{bl}{rgb}   {0.,   0.5,   1. }
\definecolor{mg}{rgb}   {0.85,  0.,    0.85}
\definecolor{yl}{rgb}   {0.8,  0.7,   0.}
\definecolor{or}{rgb}  {0.7,0.2,0.2}

\newcommand{\wh}{\widehat}

\newtheorem{theorem}{Theorem} [section]
\newtheorem{lemma}[theorem]{Lemma}
\newtheorem{proposition}[theorem]{Proposition}
\newtheorem{remark}[theorem]{Remark}

\newtheorem{definition}[theorem]{Definition}
\newtheorem{corollary}[theorem]{Corollary}

\newcommand{\noi}{\noindent}

\newcommand{\R}{\mathbb{R}}

\newcommand{\T}{\mathbb{T}}

\newcommand{\cF}{\mathcal{F}}

\let\Re=\undefined\DeclareMathOperator*{\Re}{Re}

\newcommand{\E}{\mathbb{E}}

\newcommand{\al}{\alpha}
\newcommand{\be}{\beta}
\newcommand{\dl}{\delta}

\newcommand{\Dl}{\Delta}
\newcommand{\eps}{\varepsilon}
\newcommand{\kk}{\kappa}

\newcommand{\ld}{\lambda}
\newcommand{\Ld}{\Lambda}
\newcommand{\s}{\sigma}

\newcommand{\wt}{\widetilde}

\newcommand{\dd}{\partial}

\renewcommand{\l}{\ell}
\renewcommand{\o}{\omega}
\renewcommand{\O}{\Omega}

\newcommand{\ind}{\mathbf 1}

\newcommand{\PP}{\mathbb{P}}

\renewcommand{\H}{\mathcal{H}}

\def\sgn{\textup{sgn}}

\newcommand{\ov}{\overline}

\numberwithin{equation}{section}
\numberwithin{theorem}{section}

\makeatletter
\@namedef{subjclassname@2020}{\textup{2020} Mathematics Subject
 Classification}
\makeatother

\newcommand{\cS}{\mathcal{S}}

\newcommand{\cX}{\mathcal{X}}

\newcommand{\cZ}{\mathcal{Z}}
\newcommand{\fA}{\mathfrak{A}}

\newcommand{\fS}{\mathfrak{S}}

\newcommand{\Gam}{\Gamma}

\newcommand{\ta}{\theta}
\newcommand{\AND}{\quad \text{and}\quad}

\DeclareMathOperator{\Dom}{Dom}

 \newcommand{\bfY}{\mathbf{Y}}
\newcommand{\bfX}{\mathbf{X}}
\newcommand{\sfG}{\mathsf{G}}
\newcommand{\sfK}{\mathsf{K}}

\newcommand{\rmd}{\mathrm{d}}
\newcommand{\bsy}{\boldsymbol{y}}

 \newcommand{\pmxi}[1]{\pmb{\xi_{#1}}}
\newcommand{\pmeta}[1]{\pmb{\eta_{#1}}}
\newcommand{\pmal}[1]{\pmb{\al_{#1}}}
\newcommand{\pmx}[1]{\pmb{x_{#1}}}
\newcommand{\pmy}[1]{\pmb{y_{#1}}}

\newcommand{\pms}[1]{\pmb{s_{#1}}}
\newcommand{\pmD}[1]{\pmb{D_{#1}}}
\newcommand{\pmwtD}[2][]{\pmb{\wt{D}_{#2}#1}}

\begin{document}
\baselineskip=14pt

\title[HAM:  solvability and asymptotics]
{Hyperbolic Anderson model with time-independent 
 rough noise:
   solvability and exact asymptotics}

\author[X. Wang, P. Xia, and G. Zheng]
{Xiong Wang, Panqiu Xia, and Guangqu Zheng}

\address{Xiong Wang,
School of Mathematics\\
Sun Yat-sen University\\
Xingang West Road\\
Guangzhou, Guangdong 510275\\
China}
\email{wangx937@mail.sysu.edu.cn or xiongwang@ualberta.ca}

\address{Panqiu Xia,
School of Computational and Mathematical Sciences\\
Cardiff University\\
Senghennydd Road, Cathays\\
Cardiff, Wales CF24 4AG\\
UK}
\email{xiap@cardiff.ac.uk}

\address{Guangqu Zheng,
Department of Mathematics and Statistics\\
Boston University\\
665 Commonwealth Avenue\\
Boston, MA 02215\\
USA}
\email{gzheng90@bu.edu}



\begin{abstract}
We study the one-dimensional hyperbolic Anderson model driven by time-independent Gaussian noise in the Skorokhod sense.
The noise has spectral density $(2\pi)^{-1}|\xi|^{1-2H}$ with $H<1/2$ and is therefore spatially rougher than white noise.
 For $0<H<1/2$, it coincides, up to a multiplicative constant, with the distributional derivative of a two-sided fractional Brownian motion with Hurst parameter $H$. 
We identify a necessary and sufficient condition for the existence of a unique random-field solution:  $H > - 1/4$. In particular, this range includes negative values of $H$, extending beyond the classical fractional Brownian motion framework. Throughout the solvability range, we identify the exact leading constant in the joint large-time and high-moment asymptotics of $\log\|u(t,x)\|_p$.

\end{abstract}

\date{\today}
\maketitle

\vspace{-10mm}

\tableofcontents

\section{Introduction}\label{sec_intro}

We consider the one-dimensional hyperbolic Anderson model (HAM) driven by
time-independent Gaussian noise:
\begin{equation}\label{ham}
\begin{dcases}
 \frac{\dd^2 u}{\dd t^2}(t,x)
 =\frac{\dd^2 u}{\dd x^2}(t,x)
 +u(t,x)\diamond\dot W(x),
 &(t,x)\in(0,\infty)\times\R,\\
 u(0,x)=1,\qquad
 \frac{\dd u}{\dd t}(0,x)=0,
 &x\in\R.
\end{dcases}
\end{equation}
Here $\diamond$ denotes the Wick product, indicating the stochastic term
is interpreted in the Skorokhod sense; see, e.g.,
\cite[Section~1.3]{Nua06}.
The noise $\dot W$ is centered and spatially homogeneous, with
covariance
\begin{equation}\label{cov}
 \E[W(\phi)W(\psi)]
 =\frac{1}{2\pi}\int_{\R}
 \wh\phi(\xi)\ov{\wh\psi(\xi)}
 |\xi|^{1-2H}\,d\xi,
\end{equation}
where $H<1$ and $\phi,\psi\in\cS(\R)$ are real-valued
Schwartz functions. We use the Fourier transform convention
\begin{equation}\label{four}
\cF\phi(\xi) = \wh\phi(\xi)\coloneqq \int_{\R}e^{ix\xi}\phi(x)\,dx.
\end{equation}
For $H\in(0,1)$, the noise is, up to a multiplicative constant,
the distributional derivative of a two-sided fractional Brownian
motion with Hurst parameter $H$. The spectral definition
\eqref{cov} also makes sense for $H\le0$, although $H$ then
no longer has this fractional Brownian interpretation.
Our interest lies in the regime $H<1/2$, where the noise is
spatially rougher than white noise.

We seek a random-field solution to \eqref{ham} satisfying
\begin{equation}\label{mild_w}
 u(t,x)
 =1+\int_{\R}
 \bigg(\int_0^t G_{t-s}(x-y)u(s,y)\,ds\bigg)W(dy),
\end{equation}
where the outer integral is a Skorokhod integral and $G$ is the
fundamental solution of the one-dimensional wave equation:
\[
 G_t(x)\coloneqq \frac12\ind_{\{|x|<t\}} \quad 
\text{and hence} \quad 
 \wh G_t(\xi)=\frac{\sin(t|\xi|)}{|\xi|},
 \qquad t>0.
\]

\noi
We refer readers to Section \ref{sec_prelim} for a brief review of the Gaussian framework and
the associated solution theory.

The random-field approach to stochastic partial differential
equations (SPDEs) originates in Walsh's Saint-Flour lecture notes
\cite{Wal86}, which developed stochastic integration with respect
to martingale measures and its applications to SPDEs.
Dalang \cite{Dal99} extended this framework to a broader class
of spatially homogeneous Gaussian noises that are white in time.
A central role is played by the spectral integrability condition
\begin{equation*}
 \int_{\R^d}\frac{\mu(d\xi)}{1+|\xi|^2}<\infty,
 \tag{\textit{Dalang's condition}}
\end{equation*}
where $\mu$ is the spectral measure of the   spatial covariance.
When the spatial covariance is a nonnegative distribution,
and hence a Radon measure, this condition ensures existence
and uniqueness of random-field solutions to the one-dimensional
stochastic heat and wave equations with globally Lipschitz
coefficients and zero initial data; see
\cite[Theorem 13 and Remark 14(a)]{Dal99}.
It is also necessary in the additive-noise case; see
\cite[Theorem 11 and Remark 12(a)]{Dal99}.

The classical Walsh--Dalang theory relies on the martingale
structure induced by temporal white noise and therefore does
not apply directly to temporally correlated noise.
The Skorokhod integral, defined through Malliavin calculus,
provides a natural alternative.
For general nonlinear SPDEs driven by temporally correlated
Gaussian noise, however, the corresponding solution theory
is substantially less developed.
The class of \emph{Anderson models}, in which the noise acts
linearly and multiplicatively on the solution, is more tractable.
For these models, iterating the mild equation yields a formal
series of multiple stochastic integrals, which can be made
rigorous through Wiener--It\^o chaos expansions under suitable
assumptions on the noise covariance.
This approach has been used to construct Skorokhod solutions
for both the \emph{parabolic Anderson model}
\cite{HHNT15} and the \emph{hyperbolic Anderson model}
\cite{Bal12,BS17}.

A widely studied class of driving noises is fractional Gaussian
noise, obtained by differentiating a fractional Brownian sheet
with temporal Hurst parameter $H_0$ and spatial Hurst parameters
$H_1,\ldots,H_d$.
The case $H_0=1/2$ corresponds to temporally white noise.
These parameters have a direct regularity interpretation:
coordinate sections of the underlying sheet have locally
H\"older continuous paths of every order below the corresponding
Hurst parameter.
Solvability conditions expressed in terms of these parameters
therefore describe how temporal and spatial regularities interact
with the evolution operator.
This reflects a familiar theme in PDE theory: greater regularity
in one variable may compensate for lower regularity in another,
with the precise balance depending on the equation.

For the hyperbolic Anderson model in dimensions $d\in\{1,2,3\}$,
such a trade-off was established in \cite{CDST25} under the
assumption that the spatial covariance is nonnegative.
When the temporal covariance is that of fractional Gaussian
noise with $H_0\in(1/2,1)$, the necessary and sufficient
condition for solvability, under the initial-data assumptions
of that paper, is
\[
 \int_{\R^d}
 \frac{\mu(d\xi)}
 {(1+|\xi|^2)^{H_0+\frac12}}
 <\infty.
\]
The increasing denominator exponent shows how greater temporal
regularity weakens the required spatial spectral integrability.
For fractional spatial noise with $H_j\in(1/2,1)$, this criterion
takes the form
$H_0+\sum_{j=1}^d H_j>d-\frac12.$
At the time-independent endpoint $H_0=1$, the spectral condition
remains sufficient; see \cite[Theorem~1.1]{CDST25}.
This agrees with the global solvability results of \cite{BCC22}
for the classes of time-independent spatial covariances considered
there in dimensions $d\le3$.
The critical case of three-dimensional spatial white noise behaves
differently: there exists a finite deterministic time $T_2$ such that
the Wiener chaos expansion converges in $L^2(\O)$ for $t<T_2$
but its second-moment series diverges for $t>T_2$; see
\cite[Theorem~1.7 and equation~(1.12)]{BCC22}.

The nonnegativity assumption on the spatial covariance excludes
fractional spatial noise with a Hurst parameter below $1/2$.
In this rough regime, the covariance is a distribution rather
than a Radon measure and is naturally described through spectral
measures and fractional Sobolev inner products, as in \eqref{cov}.
Dalang's condition alone need not ensure solvability of equations
with multiplicative noise; see, e.g., \cite{BJQS15,HHLNT17}.
In particular, in one spatial dimension, writing $H=H_1\in(0,1/2)$,
the condition $H>1/4$ is necessary and sufficient for
square-integrable random-field solvability of the hyperbolic
Anderson model with temporally white noise and the initial data
in \eqref{ham}; see \cite{BJQS17}.

Temporal correlations can relax this spatial regularity restriction
in the parabolic setting.
For the parabolic Anderson model with $H_0\in(1/2,1)$ and
$H\in(0,1/2)$, some sufficient conditions for solvability under suitable assumptions on the
initial data, allowing spatial Hurst parameters below $1/4$ are developed in \cite{Che19,BCM22}.
Further results on solvability with rough Gaussian noise are
developed in \cite{CDOT21,CH23,LHW24}; in particular, $2H_0+H>5/4$ is shown in \cite{LHW24} 
to be the necessary and sufficient condition when $H_0\in(1/2,1)$ and $H\in(1/20,1/2)$.

For hyperbolic Anderson models, the corresponding interaction
between temporal regularity and spatial roughness is less well
understood.
Song, Song, and Xu \cite{SSX20} established existence and
uniqueness of   Skorokhod solutions for fractional stochastic
wave equations with spatial Hurst parameter $H\in(1/4,1/2)$,
treating both temporally white noise and temporally fractional
noise with $H_0\in(1/2,1)$.
Balan and Yuan \cite{BY24} subsequently studied the time-independent
case in the same spatial Hurst range, obtaining moment estimates
and Gaussian fluctuation results.
These works accommodate temporal correlations but retain the
restriction $H>1/4$.

Related progress for equations with additive Gaussian noise
includes solvability criteria for time-fractional stochastic
diffusion equations \cite{CLX26} and sufficient and necessary conditions for the
stochastic wave equation with fractional noise \cite{LHW26}.
Such additive-noise results do not directly yield corresponding
criteria for equations driven by multiplicative noise, where one must control
the full Wiener chaos expansion.
A central question for the hyperbolic Anderson model is therefore
to what extent greater temporal regularity permits spatial
roughness beyond the range covered by existing theories.

\subsection{Main results}

In this paper, we study the time-independent hyperbolic Anderson
model in the Skorokhod setting, with particular emphasis on
solvability beyond the spatial Hurst range $H>1/4$.
In a related paper \cite{CH25}, X.~Chen and Y.~Hu  established Dalang's
condition as a necessary and sufficient criterion for solvability
in the {\it Stratonovich setting} under their hypotheses on the noise.
Under the spatial homogeneity assumption in
\cite{CH25}, they also
obtained precise moment asymptotics.
Their Theorem~1.2  requires $0<\alpha<\min\{2,d\}$ or $\alpha=d=1$.
For $d=1$ and $H<1/2$, \eqref{cov} gives $\alpha=2-2H>1$, outside this range.
Recent developments in the study of Stratonovich equations include, among others, the work of \cite{CDST23} and the pioneering contribution of L. Chen, which explores the use of AI-assisted methods in mathematical research \cite{Che26}.

Our second objective is to determine the precise moment growth
of the solution and its implications for intermittency.
Informally, intermittency means that the solution develops tall peaks
in a few small spatial regions while remaining much smaller elsewhere.
At the level of moments, it is reflected in growth rates that,
after normalization by the moment order, increase strictly with
that order.
The systematic study of this phenomenon for SPDEs was developed
by Carmona and Molchanov \cite{CM94} for the parabolic Anderson
model and subsequently by Dalang and Mueller \cite{DM09} for
the hyperbolic Anderson model.
Further developments include
\cite{CD15-1,CD15-2,BJQS17,BS19,BCM22,BNQZ22}. For further background on intermittency and the long-time asymptotic behaviour of solutions to stochastic heat equations, we refer the reader to the monograph of Khoshnevisan \cite{Kho14} and the monograph of K\"{o}nig \cite{Kon16}, which provide complementary perspectives on these topics.

\begin{theorem}\label{thm_ham}
Let $H<1/2$. The hyperbolic Anderson model \eqref{ham}
admits a unique random-field solution in the Skorokhod sense
if and only if $H>-1/4$.
\end{theorem}

To state the moment asymptotics, we write
$\|X\|_p\coloneqq (\E|X|^p)^{1/p}$ for a random variable $X$.
For deterministic functions, $\|\cdot\|_{L^2}$ denotes the
$L^2(\R)$ norm, and
$W^{1,2}(\R)
 =\{g\in L^2(\R) \colon g'\in L^2(\R)\},$
where $g'$ is the weak derivative.
We denote by $\H$ the Hilbert space associated with the noise;
see Section~\ref{ss_FGF} for more details. 

\begin{theorem}\label{thm_mom}
Let $H\in(-1/4,1/2)$, and let $u$ be the solution to \eqref{ham}.
Then $u(t,x)\in L^p(\O)$ for every $p\ge2$ and
$(t,x)\in\R_+\times\R$.
Moreover, for every $x\in\R$,
\[
 \lim_{t_{p,H}\to\infty}
 \frac{\log\|u(t,x)\|_p}{t_{p,H}}
 =\sfK_H,
\]
where
\begin{align}\label{tpHK}
 t_{p,H}\coloneqq   (p-1)^{\frac{1}{2H+1}}
 t^{\frac{2H+2}{2H+1}},
 \qquad
 \sfK_H 
 \coloneqq \frac{2H+1}{2(H+1)^{\frac{2H+2}{2H+1}}}
 \sfG_H^{\frac{H+1}{2H+1}},
\end{align}
and
\begin{equation}\label{G_H}
 \sfG_H \coloneqq \sup\Big\{  \|g^2\|_{\H}-\|g'\|_{L^2}^2 \, \colon   \,
 g\in W^{1,2}(\R) \,\,{\rm and}\,\, \|g\|_{L^2}=1  \Big\}.
\end{equation}
The limit is joint in $p\ge2$ and $t>0$, subject only to
$t_{p,H}\to + \infty$.
\end{theorem}

In particular, Theorem~\ref{thm_mom} gives the following
large-time and high-moment asymptotics.

\begin{corollary}\label{cor_mom}
Let $H\in(-1/4,1/2)$ and $x\in\R$.
\begin{enumerate}
\item For every fixed $p\ge2$,
\[
 \lim_{t\to\infty}
 \frac{\log\|u(t,x)\|_p}
 {t^{\frac{2H+2}{2H+1}}}
 =
 \sfK_H(p-1)^{\frac{1}{2H+1}}.
\]

\item For every fixed $t>0$,
\[
 \lim_{p\to\infty}
 \frac{\log\|u(t,x)\|_p}
 {(p-1)^{\frac{1}{2H+1}}}
 =
 \sfK_H t^{\frac{2H+2}{2H+1}}.
\]
\end{enumerate}
\end{corollary}

Combining our results with \cite[Theorem~1.6]{BCC22}, we see that the solvability criterion in Theorem~\ref{thm_ham} remains valid when the restriction $H<1/2$ is relaxed to $H<1$. Similarly, the moment asymptotics in Theorem~\ref{thm_mom} and Corollary~\ref{cor_mom} extend to the full range $-1/4<H<1$. Moreover, under the normalization adopted in \cite[(1.7)]{BCC22}, spatial scaling yields
\(
\mathcal M(\gamma,1)=2^{(1-H)/(1+H)}\sfG_H,
\)
showing that the leading constant in their asymptotics coincides with the one in \eqref{tpHK}.
We emphasize, however, that the extension to the rough regime $H<1/2$ does not follow directly from \cite{BCC22}. Indeed, in our regime, the spatial covariance $\gamma$ is not a positive measure and thus  fails to admit a factorization of the form $\gamma=K*K$ with $K\geq 0$, as required by \cite[Assumption~A(iii)]{BCC22}; see also \cite[Assumption~B(ii)]{GSW26}. The same obstruction applies to the related factorizations employed in \cite{CHSX15,Che16,CH25}. In Section~\ref{SEC12}, we explain how we establish the exact moment asymptotics without relying on such factorizations or the associated positivity arguments.

Moreover, our main results may extend to all spatial dimensions \(d \geq 1\). Indeed, suppose that the spectral measure has density
\(
\mu(\xi) = (2\pi)^{-d} |\xi|^{d - 2H}
\).
Our current arguments can be readily extended to establish solvability and the corresponding asymptotic results in the range \(d-1\leq H<d\), which corresponds to \(0\leq H<1\) in dimension one. However, we expect the optimal solvability range to be
\(
H > \frac{2 d - 3}{2} \vee \frac{3 d - 4}{4}
\).
Extending the argument of Lemma \ref{lem_tail}, which bridges the gap in dimension one, to higher dimensions appears to require new ideas. We therefore leave this problem for future investigation.

\subsection{Proof strategies and remarks} \label{SEC12}

Let us first briefly sketch the proof of Theorem \ref{thm_ham}.

\smallskip
\noi
\raisebox{0.15ex}{\(\bullet\)}
{\bf Solvability.}
Following the Wiener chaos approach of \cite{Hu02}, formal iteration
of the mild equation \eqref{mild_w} yields the expansion
\begin{equation}\label{chaos-u}
 u(t,x)= 1 + \sum_{n=1}^{\infty}I_n(\wt f_{n,t,x}),
\end{equation}
where $I_n$ denotes the $n$-fold multiple Wiener--It\^o integral.
We set $f_{0,t,x}=\wt f_{0,t,x}=1$, with $I_0(1)=1$.
For $n\ge1$, $\wt f_{n,t,x}$ is the symmetrization, in the
spatial variables, of
\begin{equation}\label{ker_w}
 f_{n,t,x}(y_1,\ldots,y_n)
 \coloneqq 
 \int_{0<s_1<\cdots<s_n<t}
 \prod_{j=1}^n
 G_{s_{j+1}-s_j}(y_{j+1}-y_j)
 \,ds_1\cdots ds_n,
\end{equation}
with $(s_{n+1},y_{n+1})=(t,x)$.

In Section~\ref{ss_build}, we show that every
solution of \eqref{mild_w} has $\wt f_{n,t,x}$ as its $n$th chaos
kernel. Uniqueness follows from the orthogonality relation \eqref{iso}. See
\cite[Theorem~2.2]{BCC22} for time-independent noise and
\cite[Proposition~3.2]{BJQS17} for temporally white noise.
For the covariance \eqref{cov}, Section~\ref{sec_proof} proves the
converse by constructing the solution and verifying the domain
condition in Definition~\ref{def_sol}.

\begin{proposition}\label{prop_chaos}
Let $H<1/2$ and $T\in(0,\infty]$.
Equation \eqref{ham} admits a square-integrable Skorokhod
random-field solution on $[0,T)\times\R$ if and only if
\begin{equation}\label{cond_w}
 \sum_{n=1}^{\infty}
 n!\|\wt f_{n,t,x}\|_{\H^{\otimes n}}^2<\infty
 \quad\text{for every } \,  (t,x)\in[0,T)\times\R.
\end{equation}
In that case, the solution is unique and is given by
\eqref{chaos-u}, with convergence in $L^2(\O)$
and
$\E\big[ |u(t,x)|^2 \big]
 =
1 + \sum_{n=1}^{\infty}
 n!\|\wt f_{n,t,x}\|_{\H^{\otimes n}}^2.
$
\end{proposition}

Section~\ref{sec_proof} proves both directions of
Theorem~\ref{thm_ham}. When $H\leq-1/4$, the symmetrized
second-chaos kernel has infinite squared norm for every $t>0$.
Indeed, restricting the Fourier integral to a selected region gives a lower bound
proportional to $\int_{4/t}^\infty a^{-2-4H}\,da$. The divergence of this integral rules out an
$L^2$ solution by \eqref{cond_w}.

For $-1/4<H<1/2$, Lemma~\ref{lem_wker} starts by taking the
Laplace transform in the terminal time. Writing
$\eta_j=\xi_1+\cdots+\xi_j$, we obtain one wave resolvent for
each $\eta_j$. After Plancherel in time, we bound each spatial
weight $|\eta_j-\eta_{j-1}|^{1-2H}$ by a sum of powers of
$|\eta_j|$ and $|\eta_{j-1}|$. Each $\eta_j$ occurs in at most
two adjacent weights, so the largest power is $2(1-2H)<3$.
This puts all the resulting one-dimensional integrals within
the range of the resolvent estimate \eqref{scalar_resolvent_bd}.
Time scaling then gives the factorial bound \eqref{wker_bd}.

Since full and partial symmetrization are contractions, the
same bound gives convergence of the solution's chaos series
and controls the time convolution in $L^2(\O;\H)$. Closedness
of the Skorokhod integral then justifies passing to the limit
in \eqref{mild_w}.
\medskip

\noi
\raisebox{0.15ex}{\(\bullet\)}
{\bf Moment asymptotics.}
For the lower bound, we test $u(t,x)$ against the Wick exponential
$\cZ$ associated with a real Schwartz potential $V\in\H^*$.
By \eqref{S_trans}, the pairing $\E[u(t,x)\cZ]$ equals $w_V(t,x)$,
where $w_V$ solves the deterministic wave equation \eqref{wave_pot}.
In \cite[(6.9)--(6.11)]{BCC22}, choosing $f*\gamma\geq0$ makes the
Wick pairing termwise nonnegative. However, in our setting every nonzero admissible $V$
changes sign when $H\leq0$, since \eqref{Hast} forces $\int_\R V(x)\,dx=0$.
Let $\Ld(V)$ denote the variational quantity defined in \eqref{spec_val}.
If $\Ld(V)>0$, Lemma~\ref{wave_grow} gives
$
 \sup_{z\in\R}|w_V(t,z)|
 \geq \cosh\bigl(t\sqrt{\Ld(V)}\bigr)$,
$t\geq0$.
 The proof tests the wave equation against a positive
eigenfunction of $-\partial_x^2-V$ corresponding to its lowest
eigenvalue $-\Ld(V)$.
Next, H\"older's inequality and spatial stationarity transfer this growth back to our main object 
$\|u(t,x)\|_p$ in \eqref{holder_sup}. Rescaling $V(y)\mapsto\ld^2V(\ld y)$
balances wave growth and Gaussian cost, giving the coefficient
$\sqrt{\Ld(V)}-\tfrac12\|V\|_{\H^*}^2$ of $t_{p,H}$ in \eqref{prof_bd}.
At a fixed $\|V\|_{\H^*}$, duality \eqref{dual_id} reduces the spectral
optimization to the variational problem \eqref{M_exp} for $\sfG_H$.
Optimizing over $\|V\|_{\H^*}$ as in Lemma~\ref{common_C} gives the desired lower
bound $\|u(t,x)\|_p\geq\tfrac12\exp \bigl(\sfK_H t_{p,H}\bigr)$
in \eqref{mom_bounds}.

For the upper bound, we apply hypercontractivity \eqref{chaos_hyp}
and time scaling \eqref{coeff_scale} to estimate $\|u(t,x)\|_p$ by a
weighted sum of chaos kernel norms at time one. In \cite[Lemmas~3.5 and~3.7]{BCC22},
$\gamma=K*K$ with $K\geq0$ gives nonnegative kernels that are
nondecreasing in time, allowing a reverse Cauchy--Schwarz estimate.
Our covariance admits no such factorization for $H<1/2$.
First, we use \eqref{freq_int} to bound the temporal-frequency integral by the zero-frequency value
$F_n$ from \eqref{def_Dn} times $C_H n^{\frac{3}{2}-H}$, whose $n$-th root tends to one.
Next, squaring the Laplace transform gives $F_n$ as a multiple of
$\E[\bfY^n]$ by \eqref{eq_Y-moment}, with $\bfY=Y(\tau,\tau')$ as in \eqref{def_YR}--\eqref{def_Y}
for independent Brownian local times and i.i.d.\ $\operatorname{Exp}(1)$ times.
Then, Lemma~\ref{lem_tail_Y} controls the cutoff error $\bfY-\bfY_R$.
For $H\leq0$, the single-path energy \eqref{def_XR} has mean
$\E[\bfX_R]=(4/\pi)\int_0^R \xi^{1-2H}(2+\xi^2)^{-1}\,d\xi\uparrow\infty$.
The Mukherjee--Varadhan large deviation principle for Brownian occupation
measures \cite[Theorem~4.1]{MV16} assigns $g^2$ the cost
$\tfrac12\|g'\|_{L^2}^2$. Varadhan's lemma applies to the continuous
cutoff interaction in \eqref{ieq_ldp}. Balancing $\ld\|g^2\|_\H$ against
this cost gives $\sfG_H$ in \eqref{band_lap} via \eqref{M_exp}. Proposition~\ref{prop_exp} yields
$\E[\exp(\ld\sqrt{|\bfY|})]<\infty$ for $0<\ld<(2/\sfG_H)^{(1+H)/2}$.
The $2n$-th term bounds $\E[|\bfY|^n]$, yielding \eqref{ieq_claim}.
Finally, summing by \eqref{series-asymp} gives
$\|u(t,x)\|_p\leq C_{H,\eps}e^{(\sfK_H+\eps)t_{p,H}}$ for each $\eps>0$
in \eqref{mom_bounds}.

\medskip

We conclude with three remarks. The first two compare our solvability
criterion and exact moment asymptotics with earlier results. The third
identifies the critical parameter for exponential integrability of the two-path 
Brownian functional $\sqrt{|\bfY|}$.

\begin{remark}[The solvability threshold]\rm
For the spectral measure in \eqref{cov},
\[
 \int_\R \frac{|\xi|^{1-2H}}{(1+\xi^2)^s}\,d\xi<\infty
 \quad\Longleftrightarrow\quad H>1-s,
\]
where $H<1$ and $s>0$.
Indeed, the integrand is integrable near zero and behaves like
$|\xi|^{1-2H-2s}$ at infinity. Thus Dalang's condition ($s=1$)
holds precisely when $H>0$, while Theorem~\ref{thm_ham} gives a
global Skorokhod solution also for $-1/4<H\leq0$.

For time-independent noise, \cite[Theorem~1.1]{CDST25} gives
a sufficient condition for solvability with $s=3/2$ under the
assumption that the spatial covariance $\gamma$ is nonnegative.
For the family \eqref{cov}, this assumption forces $H\geq1/2$.
The spectral integral with $s=3/2$, however, is finite for every
$H>-1/2$. Hence, for $-1/2<H\leq-1/4$, this integral is finite
even though the symmetrized second-chaos kernel has infinite
squared norm at every positive time. This shows that the spectral
condition alone is insufficient for solvability in the rough
regime.
\end{remark}
\begin{remark}[Exact moment asymptotics for rough noise]\rm
R.~M.~Balan, L.~Chen, and X.~Chen \cite[Theorem~1.6]{BCC22} obtained
joint large-time and high-moment asymptotics for homogeneous
covariances satisfying their Assumption~A and for spatial white
noise in dimensions one and two. In their notation, the spectral
measure satisfies $\mu(cA)=c^\alpha\mu(A)$ for $c>0$. For
\eqref{cov}, we have $\alpha=2-2H$, so their moment scale becomes
$ (p-1)^{1/(3-\alpha)}t^{(4-\alpha)/(3-\alpha)}=t_{p,H}.$
Theorem~\ref{thm_mom} gives the exact constant $\sfK_H$ on this
scale throughout $-1/4<H<1/2$, where $\alpha>1$, the spatial
dimension. In this range the covariance takes negative values on
some nonnegative test functions. This rules out the representation
$\gamma=K*K$ with $K\geq0$ used in
\cite[Assumption~A(iii)]{BCC22}.

The lower bound illustrates the role of signs when $H\leq0$.
By \eqref{Hast}, every admissible Schwartz potential $V$ satisfies
$\wh V(0)=\int_\R V(x)\,dx=0$, so every nonzero such $V$ changes
sign. When $\Ld(V)>0$, averaging the deterministic wave against
the positive ground state of $-d^2/dx^2-V$ gives the lower bound
in \eqref{holder_sup}. For the upper bound, \eqref{freq_int}
compares the integral over all temporal frequencies with its
zero-frequency value $F_n(0)$ at a cost polynomial in $n$.
This factor disappears upon taking $n$th roots.
\end{remark}

\begin{remark}[The critical Brownian exponential moment]\rm
The wave lower bound also shows that the parameter $\ld_\ast$ in
Proposition~\ref{prop_exp} is sharp:
\[
 \sup\bigl\{\ld>0 \colon
   \E[\exp(\ld\sqrt{|\bfY|})]<\infty\bigr\}
 =\ld_\ast
 =\Big(\frac{2}{\sfG_H}\Big)^{(1+H)/2},
\]
where $\sfG_H$ is defined as in \eqref{G_H}.
Proposition~\ref{prop_exp} gives finiteness for
$0<\ld<\ld_\ast$. If the exponential moment were finite for
some $\ld>\ld_\ast$, the exponential-series estimate in the
proof of Lemma~\ref{claim} would give
\[
 \limsup_{n\to\infty}
 \Big(\frac{\E[|\bfY|^n]}{(n!)^2}\Big)^{1/n}
 \leq \frac{4}{\ld^2}<\frac{4}{\ld_\ast^2}.
\]
Using this bound in \eqref{brown_bd} and summing the chaos series
as in Section~\ref{ss_upper} with $p=2$  yield
\[
 \limsup_{t\to\infty}
 \frac{\log\|u(t,0)\|_2}{t_{2,H}}
 \leq \sfK_H
      \Big(\frac{\ld_\ast}{\ld}\Big)^{2/(1+2H)}
 <\sfK_H,
\]
contradicting the lower bound in \eqref{mom_bounds}.
\end{remark}

\smallskip
\noi \raisebox{0.15ex}{\(\bullet\)} {\bf Organization of the paper.}
Section~\ref{sec_prelim} introduces the Gaussian framework and
establishes the analytic estimates for the chaos kernels and the
deterministic wave equation.
Section~\ref{ss_occ} proves the moment formula and the exponential
integrability bound for the Brownian occupation functional $\bfY$.
Section~\ref{sec_proof} proves the sharp solvability criterion in
Theorem~\ref{thm_ham}.
Section~\ref{sec_mom} proves the matching lower and upper moment
bounds in Theorem~\ref{thm_mom}.

\section{Preliminaries}\label{sec_prelim}  

 For $n\geq1$ and $t>0$, define the ordered simplex
\[
 \T_t^n
 \coloneqq \{(s_1,\ldots,s_n)\in(0,t)^n \colon
      0<s_1<\cdots<s_n<t\}.
\]

\noi
We write $\pmx{n}=(x_1,\ldots,x_n)\in\R^n$ and
$\pms{n}=(s_1,\ldots,s_n)\in\R_+^n$,
with $\R_+=[0,\infty)$.
We use the Fourier transform convention \eqref{four}.
With this normalization, Plancherel's identity reads
\begin{align}\label{plan}
 \int_\R|f(x)|^2\,dx
 =\frac1{2\pi}\int_\R|\wh f(\xi)|^2\,d\xi,
 \quad f\in L^2(\R),
\end{align}
where the Fourier transform is understood in the $L^2$ sense.

Throughout the paper, $C$ denotes a positive constant whose
value may change from line to line. A subscript, as in
$C_{\ta_1,\ldots,\ta_k}$, indicates the parameters on
which the constant may depend; it is independent of the other
variables under consideration.

\subsection{The Gaussian Hilbert space and covariance duality}\label{ss_FGF}
We collect the Fourier-space definitions and duality identities used below.
Throughout this section, 
let $H<1/2$ and set $\al=H- \frac12$.

Let $\cS(\R)$ denote the real-valued Schwartz space. Following the notation of \cite[Section~2.1]{LSSW16}, for
$\be\in\R$, set
\begin{equation*} 
\mathring H^\be(\R)
\coloneqq
\big\{\phi\in\cS(\R) \colon |\xi|^\be\wh\phi(\xi)\in L^2(\R)\big\},
\end{equation*}
and let $\dot H^\be(\R)$ be its completion under the inner product
\begin{equation*} 
\langle\phi,\psi\rangle_{\dot H^\be(\R)}
\coloneqq\frac1{2\pi}\int_\R
\wh\phi(\xi)\ov{\wh\psi(\xi)}|\xi|^{2\be}\,d\xi.
\end{equation*}
We use these homogeneous Sobolev spaces through their weighted Fourier
representatives; see \cite{DPV12,LSSW16}.
The Hilbert space associated with the covariance \eqref{cov} is
\[
\H\coloneqq\dot H^{-\al}(\R)=\dot H^{\frac{1}{2} - H}(\R).
\]
Since $\frac{1}{2}-H>0$, its defining test-function space is all of $\cS(\R)$.
The noise extends by continuity to an isonormal Gaussian process
$W=\{W(\phi) \colon \phi\in\H\}$, so that
$\E[W(\phi)W(\psi)]=\langle\phi,\psi\rangle_\H
$ for any $\phi,\psi\in\H$.

Under the unweighted Fourier pairing, the continuous dual of $\H$ is
identified with $\H^*\coloneqq\dot H^\al(\R)$, with norm
\begin{align}\label{Hast}
\|V\|_{\H^*}^2
=\frac1{2\pi}\int_\R|\wh V(\xi)|^2|\xi|^{2H-1}\,d\xi.
\end{align}
This is also the Cameron--Martin space of the noise. 
We retain the notation
\begin{align}\label{Hast2}
\mathring\H\coloneqq\mathring H^{-\al}(\R)=\cS(\R),
\qquad
\mathring\H^*\coloneqq\mathring H^\al(\R)
=\{V\in\cS(\R)\colon \|V\|_{\H^*}<\infty\}.
\end{align}
The space $\mathring\H^*$ is dense in $\H^*$: it contains the real Schwartz
functions whose Fourier transforms are smooth and compactly supported away
from zero, which are dense in the corresponding weighted $L^2$ space.
For $H>0$, every Schwartz function belongs to $\mathring\H^*$; when $H\leq0$,
the finiteness condition imposes additional vanishing at zero on its Fourier
transform. The potentials used in Section~\ref{ss_lower} belong to this space.

For $V\in\H^*$, define $(-\Dl)^\al V\in\H$ by the weighted Fourier
multiplier
\[
\cF[(-\Dl)^\al V](\xi)
=|\xi|^{2\al}\wh V(\xi).
\]
This map defines an isometric isomorphism from $\H^*$ onto $\H$, whose inverse is given by the Fourier multiplier $|\xi|^{-2\al}$. In particular,
\begin{equation}\label{E_H'}
\|(-\Dl)^\al V\|_\H^2
=\frac1{2\pi}\int_\R|\wh V(\xi)|^2|\xi|^{2\al}\,d\xi
=\|V\|_{\H^*}^2.
\end{equation}
For $V\in\H^*$ and $\phi\in\H$, covariance duality gives
\begin{equation}\label{id_lap-phi-psi}
\langle(-\Dl)^\al V,\phi\rangle_\H
=\frac1{2\pi}\int_\R\wh V(\xi)\ov{\wh\phi(\xi)}\,d\xi
=\langle V,\phi\rangle_{\H^*,\H}.
\end{equation}
For $V\in\mathring\H^*$ and $\phi\in\H\cap L^2(\R)$, Plancherel's identity \eqref{plan} implies that this pairing coincides with
$\int_\R V(x)\phi(x)\,dx$.

The following lemma records a standard product estimate and the
resulting properties of the variational constant.
The scaling identity is analogous to
\cite[Lemma~2.3]{BCC22}, which treats covariance measures;
see also Lemmas B.2 and B.3 of \cite{BS19} under slightly different settings.
We include a  proof to cover the distributional covariances
considered here.

\begin{lemma}\label{lem_prod}
Let $H\in(-1/4,1/2)$, and let $\sfG_H$ be defined as in \eqref{G_H}.
For every   $g\in W^{1,2}(\R)$, we have $g^2\in\H$ and
\[
\|g^2\|_\H\leq C_H\|g\|_{L^2}^{1+H}\|g'\|_{L^2}^{1-H}.
\]
Moreover, $0<\sfG_H<\infty$, and, for every $\kk>0$,
\begin{equation}\label{M_exp}
\sfG_H(\kk)
\coloneqq
\sup\Big\{\kk\|g^2\|_\H-\|g'\|_{L^2}^2 \, \colon \, g\in W^{1,2}(\R),\ \|g\|_{L^2}=1\Big\}
=\kk^{\frac{2}{1+H}}\sfG_H.
\end{equation}

\end{lemma}

\begin{proof}
Set $s=\frac{1}{2} -H\in(0,3/4)$. 
Recall that the one-dimensional Sobolev inequality gives
\begin{align}\label{g_sup}
\|g\|_\infty
\leq \sqrt{2}\,\|g\|_{L^2}^{1/2}\|g'\|_{L^2}^{1/2}
\end{align}
for $g\in W^{1,2}(\R)$.
Indeed, for   $g\in C_c^\infty(\R)$ and $x\in\R$,
\[
|g(x)|^2
=2\int_{-\infty}^x g(y)g'(y)\,dy
\leq 2\|g\|_{L^2}\|g'\|_{L^2}.
\]
This proves \eqref{g_sup} for $g\in C_c^\infty(\R)$.
Applying the same estimate to differences of smooth approximations
extends it to $W^{1,2}(\R)$ by density.
By  \eqref{g_sup},  $g$ is bounded, and thus we have
\[
\|g^2\|_{L^2}\leq\|g\|_\infty\|g\|_{L^2}
\AND
\|(g^2)'\|_{L^2}
=2\|gg'\|_{L^2}
\leq2\|g\|_\infty\|g'\|_{L^2}.
\]
As a result,  $g^2\in W^{1,2}(\R)$.
Since $0<s<1$, H\"older's inequality with conjugate exponents
$1/(1-s)$ and $1/s$, followed by Plancherel's identity \eqref{plan}, gives
\begin{align*}
\|g^2\|_\H^2
&=\frac{1}{2\pi}\int_\R
\bigl(|\wh{g^2}(\xi)|^2\bigr)^{1-s}
\bigl(|\xi|^2|\wh{g^2}(\xi)|^2\bigr)^s\,d\xi\\
&\leq
\bigg(\frac{1}{2\pi}\int_\R
|\wh{g^2}(\xi)|^2\,d\xi\bigg)^{1-s}
\bigg(\frac{1}{2\pi}\int_\R
|\xi|^2|\wh{g^2}(\xi)|^2\,d\xi\bigg)^s
=\|g^2\|_{L^2}^{2(1-s)}\|(g^2)'\|_{L^2}^{2s}.
\end{align*}
Then, it follows from \eqref{g_sup} with $s=\frac12-H$
that 

\noi
\begin{align*}
\|g^2\|_\H
\leq \|g^2\|_{L^2}^{1-s}\|(g^2)'\|_{L^2}^s
&\leq 2^s\|g\|_\infty\|g\|_{L^2}^{1-s}\|g'\|_{L^2}^s\\
&\leq C_H\|g\|_{L^2}^{1+H}\|g'\|_{L^2}^{1-H}.
\end{align*}

\noi
Restricting to $\|g\|_{L^2}=1$ and using $1-H<2$, we conclude that
\[
\sfG_H
\leq \sup\big\{C_Hz^{1-H}-z^2 \, \colon\, z\geq0 \big\}
<\infty.
\]

\smallskip

Next, we show $\sfG_H > 0$ and prove the scaling identity \eqref{M_exp}.
For $g\in W^{1,2}(\R)$ with $\|g\|_{L^2}=1$ and $\ld>0$, 
we set
$g_\ld(x) = \sqrt\ld\,g(\ld x).$
Then, we have  
\begin{equation}\label{g-scale}
\|g_\ld\|_{L^2}=1, 
\quad \|g_\ld^2\|_\H=\ld^{1-H}\|g^2\|_\H,
\AND
\|g_\ld'\|_{L^2}^2=\ld^2\|g'\|_{L^2}^2.
\end{equation}
As $\|g^2\|_\H>0$ and $1-H<2$, the quantity
$\ld^{1-H}\|g^2\|_\H-\ld^2\|g'\|_{L^2}^2$ is positive for sufficiently
small $\ld>0$. Hence $\sfG_H>0$.
Choosing  $\ld=\kk^{1/(1+H)}$ in \eqref{g-scale}
yields
\[
\kk\|g_\ld^2\|_\H-\|g_\ld'\|_{L^2}^2
=\kk^{\frac{2}{1+H}}
\big(\|g^2\|_\H-\|g'\|_{L^2}^2\big).
\]
Taking the supremum over $g\in W^{1,2}(\R)$ with $\|g\|_{L^2}=1$ proves
\eqref{M_exp}, since $g\mapsto g_\ld$ is a bijection of this class.
Hence, the proof of Lemma \ref{lem_prod} is complete. 
\qedhere

\end{proof}

\subsection{A deterministic wave equation with potential}
\label{ss_test}

We establish a deterministic growth estimate that will be used
in the proof of the sharp lower moment bound.

For a real-valued potential $V\in\cS(\R)$, we define the Schr\"odinger operator $L_V$ by

\noi
\begin{equation*} 
L_V g \coloneqq - g'' - V g,
\end{equation*}
with domain $W^{2,2}(\R)$, the space of functions
in $L^2(\R)$ whose first and second weak derivatives belong to
$L^2(\R)$.
Since $V$ is bounded, $L_V$ is self-adjoint on $L^2(\R)$;
see, e.g.,  \cite[Theorem~6.4]{Tes14}.
Denote
\begin{equation}\label{spec_val}
\Ld(V)\coloneqq
\sup\Big\{\int_\R V(x)g(x)^2\,dx-\|g'\|_{L^2}^2\,\colon\, g\in W^{1,2}(\R),\,\|g\|_{L^2}=1 \Big\}.
\end{equation}

\noi
Let $\s(L_V)$ denote the spectrum of $L_V$.
By the variational characterization of the bottom of the spectrum
\cite[Theorem~4.12, $n=1$]{Tes14},
\begin{equation}\label{spec_bot}
\inf\s(L_V)
=
\inf_{\substack{g\in W^{2,2}(\R),\ \|g\|_{L^2}=1}}
\langle L_Vg,g\rangle_{L^2(\R)}
=
-\Ld(V).
\end{equation}
The last equality follows from the integration-by-parts identity
\[
\langle L_Vg,g\rangle_{L^2(\R)}
=\|g'\|_{L^2}^2-\int_\R V(x)g(x)^2\,dx
\]
and the density of $W^{2,2}(\R)$ in $W^{1,2}(\R)$.
The identity holds first for smooth compactly supported functions
and extends to $W^{2,2}(\R)$ by approximation.

Consider the wave equation
\begin{equation}\label{wave_pot}
\begin{dcases}
\Big(\frac{\dd^2}{\dd t^2}+L_V\Big)w_V(t,x)=0,
& (t,x)\in\R_+\times\R,\\[2pt]
w_V(0,x)=1,\qquad \dd_tw_V(0,x)=0,
& x\in\R.
\end{dcases}
\end{equation}
Since the initial displacement $w_V(0,\cdot)\equiv1$ does not belong to
$L^2(\R)$, we interpret $L_Vw_V$ as the differential expression
$-\frac{\partial^2}{\partial x^2}w_V-Vw_V$.

Picard iteration yields the series representation
\begin{equation}\label{wave_ser}
w_V(t,x)
=1+\sum_{n=1}^{\infty}
\int_{\R^n}
f_{n,t,x}(\bsy_n)
\prod_{j=1}^n V(y_j)\,\rmd\bsy_n,
\end{equation}
where the kernels $f_{n,t,x}$ are defined as in \eqref{ker_w}.
Since $G_t\geq0$ and $\int_\R G_t(x)\,dx=t$, induction
on the Picard terms yields the bound
$\|V\|_\infty^n t^{2n}/(2n)!$ for the absolute value
of the $n$th term, uniformly in $x\in\R$. The induction
step uses the identity
$\int_0^t(t-s)\frac{s^{2n-2}}{(2n-2)!}\,ds
 =\frac{t^{2n}}{(2n)!}$
for $n\geq 1.$
Thus the series converges uniformly on every finite time strip
$[0,T]\times\R$ and defines the unique mild solution bounded
on such strips. In particular,
\[
\sup_{x\in\R}|w_V(t,x)|
\leq
\sum_{n=0}^{\infty}
\frac{\|V\|_\infty^n t^{2n}}{(2n)!}
=
\cosh\bigl(t\sqrt{\|V\|_\infty}\bigr).
\]
Since $V$ is smooth with bounded derivatives, differentiating
the mild equation also shows that $w_V$ is a classical solution
whose derivatives up to order two are bounded on every finite
time strip.

\begin{lemma}\label{wave_grow}
Let $V\in\cS(\R)$ be real-valued and satisfy $\Ld(V)>0$.
Then, for every $t\geq0$,
\begin{equation}\label{wave_sup}
\sup_{x\in\R}|w_V(t,x)|
\geq\cosh\bigl(t\sqrt{\Ld(V)}\bigr).
\end{equation}
\end{lemma}

\begin{proof}
Since $V\in\cS(\R)$, the essential spectrum of $L_V$ is
$[0,\infty)$; see \cite[Theorem~10.2]{Tes14}.
Together with \eqref{spec_bot} and $\Ld(V)>0$, this shows
that $-\Ld(V)$ is a discrete eigenvalue.
By \cite[Theorem~10.13]{Tes14}, it admits a strictly positive
eigenfunction $\psi\in W^{2,2}(\R)$ satisfying
\[
-L_V\psi=\Ld(V)\psi
\AND
 \|\psi\|_{L^2}=1.
\]

\noi
The eigenfunction is smooth and decays exponentially.
For the decay, we can apply \cite[Theorem~C.3.4, p.~498]{Sim82}
to $L_V/2=-\frac12\,d^2/dx^2-V/2$, which matches
the normalization used there and has $\psi$ as an eigenfunction
with discrete eigenvalue $-\Ld(V)/2$.
The identity $\psi''=(\Ld(V)-V)\psi$ gives the same decay
for $\psi''$.
Since $\psi'\in W^{1,2}(\R)$ vanishes at both infinities,
integrating $\psi''$ from the appropriate infinity also gives
exponential decay of $\psi'$.
In particular, $\psi,\psi',\psi''\in L^1(\R)$.

Define
\[
F(t)\coloneqq\int_\R\psi(x)w_V(t,x)\,dx.
\]

\noi
The bounds established above for $w_V$ and its derivatives
on finite time strips, together with
$\psi,\psi',\psi''\in L^1(\R)$ and the decay of $\psi,\psi'$,
justify differentiating twice under the integral and
integrating twice by parts, with vanishing boundary terms.
Hence
\begin{align*}
F''(t)
&=-\int_\R\psi(x)L_Vw_V(t,x)\,dx
=-\int_\R L_V\psi(x)w_V(t,x)\,dx
=\Ld(V)F(t).
\end{align*}
The initial conditions in \eqref{wave_pot} imply 
$F(0)=\|\psi\|_{L^1(\R)}$, $F'(0)=0$,
and therefore we have
\begin{align*}
F(t)=\|\psi\|_{L^1(\R)}\cosh\Big(  t\sqrt{\Ld(V)}\Big).
\end{align*}
Since $\psi>0$, we conclude that
$\|\psi\|_{L^1(\R)}\cosh\bigl(t\sqrt{\Ld(V)}\bigr)
=F(t)
\leq\|\psi\|_{L^1(\R)}\sup_{x\in\R}|w_V(t,x)|.
$
Dividing by $\|\psi\|_{L^1(\R)}>0$ proves \eqref{wave_sup}.
\end{proof}

\subsection{Wiener chaos and Skorokhod solutions}\label{ss_chaos_sol}

We recall the Wiener chaos identities used to construct the solution
and estimate its moments. Write $\H^{\odot n}$ for the symmetric
tensor product of $n$ copies of $\H$. For $n\geq1$, let $I_n$ denote
the multiple Wiener--It\^o integral of order $n$, satisfying the 
following orthogonality relation

\noi
\begin{align}\label{iso}
 \E[I_m(f)I_n(g)]
 =\begin{cases}
 n!\langle f,g\rangle_{\H^{\otimes n}},&m=n,\\
 0,&m\ne n,
 \end{cases}
\end{align}
where $f\in\H^{\odot m}$ and $g\in\H^{\odot n}$.
For $f\in\H^{\otimes n}$, we set $I_n(f)=I_n(\wt f)$.
It is well known that 
every $F\in L^2(\O,\s\{W\},\PP)$ admits a unique $L^2$-expansion
\begin{align*}
 F=\E[F]+\sum_{n=1}^{\infty}I_n(f_n),
 \quad f_n\in\H^{\odot n}
\end{align*}
with
\[
 \E[F^2]= (\E[F])^2+\sum_{n=1}^{\infty}n!
 \|f_n\|_{\H^{\otimes n}}^2;
\]
see \cite[Theorem~1.1.1 and (1.9)]{Nua06}.
By hypercontractivity \cite[Corollary~2.8.14]{NP12},  
we have for any finite
$p\geq2$ and $f\in\H^{\odot n}$,
\begin{align}\label{chaos_hyp}
 \|I_n(f)\|_p
 \leq(p-1)^{n/2}\|I_n(f)\|_2
 =(p-1)^{n/2}\sqrt{n!}\,\|f\|_{\H^{\otimes n}}.
\end{align}

The Skorokhod integral $\dl$ is the adjoint of the Malliavin
derivative and is a densely defined closed operator from $L^2(\O;\H)$ to
$L^2(\O)$. If $v\in L^2(\O;\H)$ has the chaos expansion
$v=\sum_{n=0}^{\infty}I_n(v_n)$, where
$v_n\in\H^{\odot n}\otimes\H$, then $v\in\Dom\dl$ if and
only if the series
\begin{align}\label{div}
 \dl(v)\coloneqq \sum_{n=0}^{\infty}I_{n+1}(\wt v_n)
\end{align}
converges in $L^2(\O)$.
 Here $\wt v_n$ denotes symmetrization
in all $n+1$ variables. We also write
$\int_\R v(y)W(dy)=\dl(v)$; see
\cite[Definition~1.3.1, Proposition~1.3.7]{Nua06}.

\begin{definition}[Skorokhod random-field solution]\label{def_sol}
A jointly measurable, $\s\{W\}$-measurable random field $u$ is a
\emph{random-field solution} to \eqref{ham} if it is locally
$L^2(\O)$ bounded on finite space--time rectangles and satisfies
the following conditions.
\begin{enumerate}
\item  {\rm\textsf{(Initial conditions)}}
For every $x\in\R$,
\[
 u(0,x)=1,
 \qquad \lim_{t\downarrow0}\|u(t,x)-1\|_2=0.
\]
For every $\varphi\in\cS(\R)$ and $t>0$, the spatial pairing
$\int_\R\varphi(x)(u(t,x)-1)\,dx$ exists as an
$L^2(\O)$-valued Bochner integral, and
\[
 \lim_{t\downarrow0}
 \bigg\|\frac1t\int_\R\varphi(x)(u(t,x)-1)\,dx\bigg\|_2=0.
\]

\item   {\rm\textsf{(Time convolution)}}
For every $t>0$ and $x\in\R$, form the physical convolution
\[
 \Psi_{t,x}(y)=\int_0^tG_{t-s}(x-y)u(s,y)\,ds,
 \qquad y\in\R,
\]
as an $L^2(\O)$-valued Bochner integral in time.
Its Fourier transform in $y$ represents an element of
$L^2(\O;\H)$, also denoted by $\Psi_{t,x}$, which belongs to
$\Dom\dl$.

\item   {\rm\textsf{(Mild equation)}}
For every $t>0$ and $x\in\R$, equation \eqref{mild_w} holds
in $L^2(\O)$.
\end{enumerate}
Uniqueness means equality in $L^2(\O)$ at each fixed $(t,x)$
within this class. A local solution on $[0,T)$ satisfies these
conditions on each compact subinterval.
\end{definition}

\begin{remark}\label{rem_order}
\rm
(1)
Condition~(i) prescribes the initial displacement in $L^2(\O)$
and the zero initial velocity in a spatially weak sense, by
testing against Schwartz functions. The latter formulation is
essential when $-1/4 < H\leq0$, since the pointwise limit
$(u(t,x)-1)/t\to0$ in $L^2(\O)$ generally fails;
see Remark~\ref{rem_init}.
Conditions~(ii)--(iii) ensure that the Skorokhod integral
$\dl(\Psi_{t,x})$ is well defined and that
$u(t,x)=1+\dl(\Psi_{t,x})$ in $L^2(\O)$.

\smallskip
\noi
(2)
The order of integration in condition~(ii) is essential:
the time convolution must be formed before applying the
spatial divergence. Indeed, when $H\leq0$, for every $r>0$,
\[
 \|G_r\|_\H^2
 =\frac1{2\pi}\int_\R
 \frac{\sin^2(r\xi)}{|\xi|^{1+2H}}\,d\xi
 =\infty.
\]
Taking expectations in \eqref{mild_w} gives $\E[u(s,y)]=1$.
For each $s<t$, the mean of $G_{t-s}(x-\cdot)u(s,\cdot)$ is
$G_{t-s}(x-\cdot)$, which has infinite $\H$-norm. 
Condition~(ii) therefore requires the time integral $\Psi_{t,x}$
to belong to $\Dom\dl$.
\end{remark}

For later estimates, we record the Fourier representation
\begin{align}\label{four_ker}
 g_{n,t,x}(\pmxi{n})
 &\coloneqq\cF f_{n,t,x}(\pmxi{n})
 =e^{ix\sum_{k=1}^n\xi_k}
 \int_{\T_t^n}\prod_{j=1}^n
 \wh G_{s_{j+1}-s_j}\Big(\sum_{i=1}^j\xi_i\Big)
 \,d\pms{n},
\end{align}
where $s_{n+1}=t$; see also \cite[(2.5)]{BCC22}.
Note that  for $\eta\in\R$ and $\Re z>0$, the time
Laplace transform of the wave kernel satisfies
\begin{align}\label{Lapa}
 \int_0^\infty e^{-zr}\wh G_r(\eta)\,dr
 =\int_0^\infty e^{-zr}\frac{\sin(r\eta)}{\eta}\,dr
 =\frac{1}{z^2+\eta^2}.
\end{align}
For real $z>0$, this calculation appears in
\cite[p.~2154, immediately after (3.6)]{CH25}.
The same identity holds for $\Re z>0$ by writing the sine as a
difference of exponentials.

\subsection{Analytic estimates}

We record the contraction property of symmetrization
and the estimates for Gamma functions and moment series used below.
Since symmetrization is an orthogonal projection,

\noi
\begin{align}\label{sym_bd}
 \|\wt f\|_{\H^{\otimes n}}\leq\|f\|_{\H^{\otimes n}},
 \quad f\in\H^{\otimes n}.
\end{align}

The following estimates control convergence of the chaos series
and the sharp exponential rate of its moment bounds.
We write $A\sim B$ when $A/B\to1$ in the indicated limit.

\begin{lemma}\label{lem_series}
Let $a>0$.
\begin{enumerate}
\item As $n\to\infty$,
\begin{align}\label{stir}
 \Gam(an+1)
 \sim(2\pi)^{\frac{1-a}{2}}\sqrt a\,
 (n!)^a a^{an}n^{\frac{1-a}{2}}.
\end{align}
Consequently, there exists $c_a>0$ such that
\begin{align}\label{gam_low}
 \Gam(an+1)\geq c_a^n(n!)^a,
 \quad n\geq1.
\end{align}

\item As $x\to\infty$,
\begin{align}\label{series-asymp}
 \sum_{n=0}^{\infty}\frac{x^n}{(n!)^a}
 \sim\frac{(2\pi)^{\frac{1-a}{2}}}{\sqrt a}\,
 x^{\frac{1-a}{2a}}\exp\big(a x^{1/a}\big).
\end{align}
In particular, there exists $C_a>0$ such that
\begin{align}\label{series}
 \sum_{n=0}^{\infty}\frac{x^n}{(n!)^a}
 \leq C_a\exp\big(C_a x^{1/a}\big),
 \qquad x\geq0.
\end{align}
\end{enumerate}
\end{lemma}

\begin{proof}
Relation \eqref{stir} follows directly from
Stirling's formula.
The asymptotic relation \eqref{series-asymp} is the asymptotic formula for a
Hadamard--Mittag-Leffler series; see
\cite[Theorem 1]{Ger12}.
\end{proof}

\section{Estimates for Brownian local times}
\label{ss_occ}

In this section, we establish the Brownian local-time estimates used
in the sharp upper moment bound of Section~\ref{ss_upper}.
Proposition~\ref{prop_exp} gives an exponential moment estimate
whose dependence on $\sfG_H$ determines the constant in that bound.
The identity \eqref{eq_Y-moment} relates the moments of the
Brownian functional to integrals used to estimate the chaos
kernels. Lemma~\ref{rmk_Y-moment} also gives two-sided bounds
for these Brownian moments.

Let $B=\{B_t\colon t\in\R_+\}$ be a standard one-dimensional Brownian motion starting from $0$, 
and let $L=\{L_t(x)\colon (t,x)\in\R_+\times\R\}$ denote its local time, formally defined by

\noi
\begin{align*}  
L_t(x)\coloneqq \int_0^t  \dl(B_s-x) ds,
\end{align*}

\noi
where $\dl$ denotes the Dirac delta distribution. 
The Fourier transform of $L_t$ admits the representation
\begin{align}\label{def_fLt}
\wh L_t(\xi)
=\int_{\R} L_t(x)e^{i\xi x}  dx
=\int_0^t e^{i\xi B_s} ds .
\end{align}

Let $B'$ be an independent copy of $B$, and let $L'$ denote
its local time. Let $\tau,\tau'$ be independent exponential
random variables with mean $1$, independent of $B$ and $B'$.
Fix $-1/4<H<1/2$. For $R>0$, set
\begin{align}\label{def_YR}
 Y_R(s,t)
 &\coloneqq\frac1{2\pi}\int_{-R}^R
 |\xi|^{1-2H}\wh L_s(\xi)
 \ov{\wh L'_t(\xi)}\,d\xi  
\AND
 \bfY_R\coloneqq Y_R(\tau,\tau').
\end{align}
We interpret the full-frequency integral as
\[
 Y(s,t)
 \coloneqq\lim_{R\to\infty}Y_R(s,t)
 \quad\text{in }L^2(\O),
 \quad s,t\geq0,
\]
and define
\begin{align}\label{def_Y}
 \bfY\coloneqq\lim_{R\to\infty}\bfY_R
 \quad\text{in }L^2(\O).
\end{align}
The existence of these limits is established in
Lemma~\ref{lem_tail_Y}. Its proof also shows that, for a jointly
measurable version of $Y$, we have $\bfY=Y(\tau,\tau')$ almost surely.
By conjugate symmetry of $\wh L_s$ and $\wh L'_t$, the
integral in \eqref{def_YR} is real-valued over any bounded symmetric
frequency set, and the limits $Y(s,t)$ and $\bfY$ are also real-valued.

\begin{remark}\label{rmk_loc-time}
\rm
The Fourier representation \eqref{def_fLt} gives
\begin{align}\label{def_fLt-bdd}
 |\wh L_t(\xi)|\leq t,
 \quad  \forall (t,\xi)\in\R_+\times\R.
\end{align}
In particular, $|Y_R(s,t)|\leq C_HR^{2-2H}st$, so the
truncations are well defined and $\bfY_R$ has moments of every order.
When $H\leq0$, the individual local times $L_t$ need not belong to $\H$;
the full-frequency integral is therefore interpreted through the
limits above, rather than as an ordinary $\H$ inner product.
\end{remark}

We first state the exponential integrability of
$\bfY=Y(\tau,\tau')$.

\begin{proposition}\label{prop_exp}
Let $-1/4<H<1/2$, and let $\sfG_H$ and $\bfY$
be defined as in  \eqref{G_H} and \eqref{def_Y}, respectively.
Set
\begin{align}\label{ieq_exp}
 \ld_\ast
 \coloneqq
 2^{\frac{1+H}{2}}\sfG_H^{-\frac{1+H}{2}}.
\end{align}
Then
$\E\big[\exp(\ld\sqrt{|\bfY|})\big]<\infty$
for every $\ld\in(0,\ld_\ast)$.

\end{proposition}

The proof follows the established large-deviation approach to
exponential asymptotics of Brownian interaction functionals.
A closely related variational formula for square-root
self-interaction energies appears in
\cite[Theorem~1.1]{Chen23SILT} and is used in
\cite[pp.~2174--2175]{CH25}.
The motivation here is that the Brownian occupation-measure
rate function assigns the cost $\frac12\|g'\|_{L^2}^2$ to a
probability density $g^2$ (when $\|g\|_{L^2}=1$).
Combined with the interaction term, this leads to the
variational bound
\[
 \sup
 \Big\{\ta\|g^2\|_\H-\frac12\|g'\|_{L^2}^2 \,\colon \,   g\in W^{1,2}(\R) \text{\, and\,\,} \|g\|_{L^2}=1\Big\}
 =\frac12(2\ta)^{\frac{2}{1+H}}\sfG_H,
 \quad \ta>0,
\]
by \eqref{M_exp}, explaining the appearance of $\sfG_H$
in the exponential-integrability threshold.

The frequency truncation in \eqref{def_YR} yields bounded
continuous spatial kernels and
makes the comparison with single-path occupation energies
legitimate; the corresponding full-frequency comparison is
unavailable when $H\leq0$.
For each fixed cutoff, we apply the occupation-measure LDP
on the compactification of translation orbits developed in
\cite[Corollary~3.3 and Theorem~4.1]{MV16}.
We then remove the cutoff by estimating the high-frequency
remainder of the two-path functional directly, using
Lemma~\ref{lem_tail_Y}.

By Cauchy--Schwarz
and $2|ab|\leq a^2 + b^2$,
we have

\noi
\begin{align}\label{ineq_YXR}
 \sqrt{|Y_R(s,t)|}
 \leq \big(X_R(s)X'_R(t)\big)^{1/4}
 \leq \frac{1}{2}
 \Big(   \sqrt{X_R(s)}+\sqrt{X'_R(t)}   \, \Big),
\end{align}
where
\begin{align*}
 X_R(s)
=
 \frac{1}{2\pi}\int_{-R}^{R}
 |\xi|^{1-2H}|\wh L_s(\xi)|^2\,d\xi
 \quad \text{and} \quad
 X'_R(t)
=
 \frac{1}{2\pi}\int_{-R}^{R}
 |\xi|^{1-2H}|\wh L'_t(\xi)|^2\,d\xi.
\end{align*}
These quantities are finite by \eqref{def_fLt-bdd}
and $H < 1/2$.
We also write
\begin{align}\label{def_XR}
 \bfX_R\coloneqq X_R(\tau)
\quad \text{and} \quad
 \bfX'_R\coloneqq X'_R(\tau'),
\end{align}
where $\tau$ and $\tau'$ are the exponential times introduced above.

We next state two auxiliary lemmas, whose proofs are deferred
until after the proof of Proposition~\ref{prop_exp}.
The first  (Lemma \ref{lem_tail_Y})  controls the high-frequency remainder, 
while the second  (Lemma \ref{lem_band}) provides the exponential estimate at a fixed frequency
cutoff.

\begin{lemma}\label{lem_tail_Y}
Assume that $-1/4<H<1/2$. The limits defining $Y(s,t)$ and
$\bfY$ exist in $L^2(\O)$ for every $s,t\geq0$.
Moreover, $Y$ admits a jointly measurable version for which
$\bfY=Y(\tau,\tau')$ almost surely.
There exists a constant $C>0$, depending only on $H$, such that
for every positive even integer $n$ and every $R\geq1$,
\begin{align}\label{ieq_tail_Y}
 \E\big[|\bfY-\bfY_R|^n\big]
 \leq C^n(n!)^2R^{-(2H+\frac12)n}.
\end{align}
Consequently, there exist constants $C_1,C_2>0$, depending only
on $H$, such that for every $R\geq1$,
\begin{align}\label{exp33}
 \E\Big[
 \exp\Big(C_1R^{H+\frac14}
 \sqrt{|\bfY-\bfY_R|}\Big)
 \Big]\leq C_2.
\end{align}

\end{lemma}

\begin{lemma}\label{lem_band}
Assume that $-1/4<H<1/2$, and let $\sfG_H$ be defined
as in  \eqref{G_H}. Then, for every $R>0$ and $\ld>0$,
\begin{align}\label{band_lap}
 \limsup_{t\to\infty}\frac1t
 \log\E\Big[\exp\big(\ld\sqrt{X_R(t)} \, \big)\Big]
 \leq
 2^{\frac{1-H}{1+H}}\ld^{\frac{2}{1+H}}\sfG_H.
\end{align}
Consequently, with $\bfX_R$ defined by \eqref{def_XR},
we have 

\noi
\begin{align}\label{band_kill}
 \E\big[\exp(\ld\sqrt{\bfX_R})\big]<\infty
 \quad\text{for every }\ld\in(0, \tfrac{1}{2}\ld_\ast),
\end{align}
where $\ld_\ast$ is defined as in  \eqref{ieq_exp}.
\end{lemma}

Lemma~\ref{lem_band} is an application of the occupation-measure
large deviation principle of Mukherjee and Varadhan
\cite[Corollary~3.3 and Theorem~4.1]{MV16} at a fixed frequency cutoff.
Closely related square-root exponential asymptotics appear in Chen
\cite[Theorem~1.1]{Chen23SILT} and Chen and Hu
\cite[pp.~2174--2175]{CH25}, under the assumption that the spatial
covariance is nonnegative. The truncated spatial kernel used here is
positive definite but need not be pointwise nonnegative; the
occupation-measure argument applies without this additional assumption.

Combining these two lemmas, we now prove
Proposition~\ref{prop_exp}.

\begin{proof}[Proof of Proposition~\ref{prop_exp}] 

Fix $\ld\in(0, \ld_\ast)$. Then, we can choose $p>1$   such that
$
p \ld/2 < \ld_\ast/2
$
and we put $p' =p/(p-1)$. 
 By   \eqref{ineq_YXR}, we can write
\[
\sqrt{|\bfY|}
\leq \sqrt{|\bfY_R|}
+ \sqrt{|\bfY-\bfY_R|}
\leq \frac12\Big(
\sqrt{\bfX_R} + \sqrt{\bfX_R'}
\Big) + \sqrt{|\bfY-\bfY_R|}.
\]
By H\"older's inequality and the independence of the two Brownian motions and their corresponding exponential times, we obtain
\begin{align*}
\E \big[\exp \big(\ld\sqrt{|\bfY|}\big)\big]
&\leq J_1 \times J_2,
\end{align*}
where
\[
J_1\coloneqq
\Big(\E\big[\exp(\tfrac{p\ld}{2}\sqrt{\bfX_R})\big]\Big)^{2/p}
\AND
J_2\coloneqq
\Big(\E\big[\exp(p'\ld\sqrt{|\bfY-\bfY_R|})\big]\Big)^{1/p'}.
\]
By Lemma~\ref{lem_band}, $J_1 < \infty$ for every $R>0$. Furthermore, by Lemma~\ref{lem_tail_Y}, we can choose $R\geq1$ sufficiently large such that
$
p' \ld \leq C_1R^{H+\frac14}
$,
which ensures that
$ 
J_2
\leq C_2^{1/p'} < \infty
$.
Combining these estimates completes the proof. 
\end{proof}

Now let us first prove  Lemma \ref{lem_tail_Y}.

\begin{proof}[Proof of Lemma~\ref{lem_tail_Y}]
 Fix $R\geq1$, $S>R$, and a
positive even integer $n$. Consider the difference
\[
 \bfY_S-\bfY_R
 =\frac1{2\pi}\int_{R<|\xi|\leq S}
 |\xi|^{1-2H}  \wh L_\tau(\xi)
 \ov{\wh L'_{\tau'}(\xi)}\,d\xi\,,
\]
which is real-valued by the conjugate symmetry of the integrand
on the symmetric set $\{R<|\xi|\leq S\}$.

For every positive integer $n$ and $\xi_1,\ldots,\xi_n\in\R$,
the exponential-time calculation in
\cite[pp.~609--610]{Che07}, with Brownian
characteristic exponent $\psi(\xi)=\xi^2/2$, gives
\begin{align*} 
 \E\bigg[\prod_{j=1}^n\wh L_\tau(\xi_j)\bigg]
 &=\sum_{\s\in\fS_n}
   \prod_{k=1}^n
   \frac{2}{2+\big|\sum_{j=1}^k\xi_{\s(j)}\big|^2}.
\end{align*}

Define
\[
 \pmD{n}(R)
 =\big\{\pmxi{n}\in\R^n\colon 
 |\xi_j|\geq R,\ 1\leq j\leq n\big\}.
\]
By independence of $(B,\tau)$ and $(B',\tau')$, and since $n$
is even,
we can get 

\noi
\begin{align}\label{ieq_Y-YR-n}
 \E\big[|\bfY_S-\bfY_R|^n\big]
 &=(2\pi)^{-n}
   \int_{\pmD{n}(R)\cap[-S,S]^n}
   \prod_{j=1}^n|\xi_j|^{1-2H} 
   \bigg(
   \sum_{\s\in\fS_n}
   \prod_{k=1}^n
   \frac{2}{2+\big|\sum_{j=1}^k\xi_{\s(j)}\big|^2}
   \bigg)^2\,d\pmxi{n}\nonumber\\
 &\leq(2\pi)^{-n}(n!)^2
   \int_{\pmD{n}(R)}
   \prod_{k=1}^n
   \bigg[
   |\xi_k|^{1-2H} 
   \bigg(\frac{2}{2+\big|\sum_{j=1}^k\xi_j\big|^2}\bigg)^2
   \bigg]\,d\pmxi{n}.
\end{align}
Here we used Cauchy--Schwarz in the permutation sum and the
permutation invariance of the integration domain and weight.

Make the change of variables
$\eta_k=\xi_1+\cdots+\xi_k$, with $\eta_0=0$, and set
\[
 \pmwtD{n}(R)
 =\big\{\pmeta{n}\in\R^n\colon 
 |\eta_k-\eta_{k-1}|\geq R,\ 1\leq k\leq n\big\}.
\]
Expanding the bound
$|\eta_k-\eta_{k-1}|^{1-2H} 
 \leq C_H(|\eta_k|^{1-2H} +|\eta_{k-1}|^{1-2H} )$
gives
\begin{align}\label{ieq_Y-YR}
 \E\big[|\bfY_S-\bfY_R|^n\big]
 \leq C_H^n(n!)^2
 \sum_{\pmal{n}\in\fA_n}J_{\pmal{n}}(R),
\end{align}
where
\begin{align}\label{def_A-n}
 \fA_n
 =\bigg\{\pmal{n}\in\{0, 1-2H,2(1-2H) \}^n\colon 
 \sum_{j=1}^n\al_j=n(1-2H)\bigg\},
\end{align}
and
\[
 J_{\pmal{n}}(R)
 =\int_{\pmwtD{n}(R)}
   \prod_{k=1}^n q_{\al_k}(\eta_k)\,d\pmeta{n},
 \qquad
 q_\al(x)=|x|^\al\bigg(\frac{2}{2+x^2}\bigg)^2.
\]

Since $0\leq\al\leq2(1-2H) <3$
for every
$\al\in\{0,1-2H ,2(1-2H)\}$,
we have 
\[
 \int_\R q_\al(x)\,dx\leq C_H 
\AND
 \int_{|x|\geq R/2}q_\al(x)\,dx
 \leq C_H R^{\al-3}.
\]
Moreover, $|x-y|\geq R$ implies that either $|x|\geq R/2$
or $|y|\geq R/2$. Thus, for
$\al,\be\in\{0, 1-2H, 2(1-2H)\}$,
nonnegativity and the preceding integral bounds give
\begin{align*}
 \int_{|x-y|\geq R}q_\al(x)q_\be(y)\,dx\,dy
 &\leq\bigg(\int_{|x|\geq R/2}q_\al(x)\,dx\bigg)
       \bigg(\int_\R q_\be(y)\,dy\bigg)\\
 &\quad+\bigg(\int_\R q_\al(x)\,dx\bigg)
       \bigg(\int_{|y|\geq R/2}q_\be(y)\,dy\bigg)\\
 &\leq C_H\big(R^{\al-3}+R^{\be-3}\big)
 \leq C_H R^{\al+\be-3}\,.
\end{align*}

Because $n$ is even, we may retain only the constraints
$|\eta_{2j}-\eta_{2j-1}|\geq R$, $1\leq j\leq n/2$,
in the definition of $\pmwtD{n}(R)$. The resulting integral
factors over these pairs, giving
\begin{align}\label{ieq_J>0}
 J_{\pmal{n}}(R)
 &\leq\prod_{j=1}^{n/2}
   \int_{|x-y|\geq R}
   q_{\al_{2j-1}}(x)q_{\al_{2j}}(y)\,dx\,dy\nonumber\\
 &\leq C_H^nR^{(\sum_{k=1}^n\al_k)- \frac{3}{2}n}
 =C_H^nR^{-(2H+\frac12)n}.
\end{align}

\noi
Combining \eqref{ieq_J>0} with \eqref{ieq_Y-YR} and
$|\fA_n|\leq3^n$, we obtain
\[
 \E\big[|\bfY_S-\bfY_R|^n\big]
 \leq C_H^n(n!)^2R^{-(2H+\frac12)n},
\]
uniformly in $S>R$.

For $n=2$, this estimate shows that $(\bfY_R)_{R\geq1}$
is Cauchy in $L^2(\O)$; we define $\bfY$ as its limit.
Passing to an almost surely convergent subsequence and applying
Fatou's lemma proves \eqref{ieq_tail_Y} for every positive even
integer $n$.

To justify the definition at fixed times, set
\[
 M_{R,S}(s,t)
 =\E\big[|Y_S(s,t)-Y_R(s,t)|^2\big].
\]
Expanding the square in $M_{R,S}(s,t)$ and using the Gaussian
characteristic function gives a time integral over
$[0,s]^2\times[0,t]^2$ with a nonnegative integrand. Thus $M_{R,S}$ is
nondecreasing in each time variable. Independence of the
exponential times gives
\begin{align*}
 \E\big[|\bfY_S-\bfY_R|^2\big]
 &=\int_0^\infty\int_0^\infty
 e^{-u-v}M_{R,S}(u,v)\,du\,dv\\
 &\geq e^{-s-t}M_{R,S}(s,t).
\end{align*}
Consequently,
$ M_{R,S}(s,t)\leq C_H e^{s+t}R^{-(4H+1)}$,
$\forall S>R\geq1.$
Hence $Y_R(s,t)$ converges in $L^2(\O)$ to a limit $Y(s,t)$,
locally uniformly in $(s,t)$. Taking the pointwise limit along
a dyadic sequence of cutoffs (with $R=2^m$, $m\geq 1$) gives a jointly measurable version:
the displayed bound makes the $L^1(\O)$ norms of its successive
increments summable at each fixed $(s,t)$.
Conditioning on the independent exponential times then shows
that $Y_R(\tau,\tau')$ converges to $Y(\tau,\tau')$ in probability.
As it also converges to $\bfY$ in $L^2(\O)$, the two limits
agree almost surely.

Finally, let $Z_R=R^{H+\frac14}\sqrt{|\bfY-\bfY_R|}$.
For every integer $k\geq1$, Cauchy--Schwarz gives
$\E[|\bfY-\bfY_R|^k]\leq\big(\E|\bfY-\bfY_R|^{2k}\big)^{1/2}$.
Applying \eqref{ieq_tail_Y} with $n=2k$ then yields
$\E[Z_R^{2k}]\leq C_H^k(2k)!$.
Therefore, for sufficiently small $c>0$, we have 
\[
 \E[e^{cZ_R}]
 \leq2\E[\cosh(cZ_R)]
 \leq2\sum_{k=0}^\infty(C_Hc^2)^k<\infty,
\]
uniformly in $R\geq1$. This proves the exponential estimate \eqref{exp33}
and completes the proof.
\end{proof}

We now turn to the proof of Lemma~\ref{lem_band}.

\begin{proof}[Proof of Lemma~\ref{lem_band}]
For a finite nonnegative Borel measure $\eta$ on $\R$, let
$[\eta]=\{\eta(\cdot-y)\colon y\in\R\}$ denote its translation orbit.
We use the compact space
\[
 \wt{\cX}
 =\Big\{\zeta=\{[\eta_j]\}_{j\in J}\colon 
 \eta_j\geq0,\ \eta_j(\R)>0,\
 \sum_{j\in J}\eta_j(\R)\leq1\Big\},
\]
introduced in \cite[Equation~(16) and Theorem~3.2]{MV16},
equipped with the metric defined in \cite[Equation~(17)]{MV16}.
Here $J$ is at most countable, possibly empty;
the order of the orbits does not matter, but the same orbit
may appear more than once.

By \cite[Theorem~4.1]{MV16}, the translation orbits of the
normalized occupation measures $t^{-1}L_t(x)\,dx$ satisfy
a large deviation principle on $\wt{\cX}$ with speed $t$.
The rate function is
\[
 \wt I(\zeta)=\frac12\sum_{j\in J}\|g_j'\|_{L^2}^2
\]
when $\eta_j(dx)=g_j^2(x)\,dx$ with
$g_j\in W^{1,2}(\R)$ for every $j$, and is infinite otherwise.

Define
\[
 K_R(x)=\frac1{2\pi}\int_{-R}^R
 e^{i\xi x}|\xi|^{1-2H}\,d\xi.
\]
Since the spectral density
$\mathbf 1_{[-R,R]}(\xi)|\xi|^{1-2H}$
is nonnegative, even, and integrable, $K_R$ is real-valued
and positive definite. Moreover, it is continuous and vanishes
at infinity by the Riemann--Lebesgue lemma. Set
\[
 \cF_R(\zeta)
 =\bigg(
 \sum_{j\in J}\int_{\R^2}
 K_R(x-y)\,\eta_j(dx)\eta_j(dy)
 \bigg)^{1/2}.
\]

Recall that the class $\mathcal F_2$ in
\cite[Section~2.2]{MV16} consists of continuous,
translation-invariant functions of two variables that vanish
as the distance between the variables tends to infinity.
Thus $(x,y)\mapsto K_R(x-y)$ belongs to $\mathcal F_2$.
Corollary~3.3 therein states that the corresponding interaction
functional extends continuously from single orbits to
$\wt{\cX}$ by summing the interactions within each orbit.
Applied to $K_R$, this gives continuity of $\cF_R^2$.
Moreover, positive definiteness of $K_R$ and the mass constraint give
\[
 0\leq\cF_R(\zeta)^2
 \leq\|K_R\|_\infty\sum_j\eta_j(\R)^2
 \leq\|K_R\|_\infty,
\]
so $\cF_R$ is bounded and continuous for each fixed $R$. 
The configuration $\{[t^{-1}L_t(x)\,dx]\}$ contains just one
orbit, so the definition of $\cF_R$ and the Fourier representation of $K_R$ give
\begin{align}\label{FRXR}
 \cF_R\big(\{[t^{-1}L_t(x)\,dx]\}\big)
 &=\bigg(\frac1{2\pi t^2}\int_{-R}^R
 |\xi|^{1-2H}|\wh L_t(\xi)|^2\,d\xi\bigg)^{1/2}
 =t^{-1}\sqrt{X_R(t)}.
\end{align}

\noi
Applying Varadhan's integral lemma
\cite[Theorem~4.3.1]{DZ10} to the occupation-measure LDP
with speed $t$ and rate function $\wt I$, and using
\eqref{FRXR},
we can obtain
\begin{align}\label{ieq_ldp}
 \lim_{t\to\infty}\frac1t
 \log\E\big[\exp(\ld\sqrt{X_R(t)})\big]
 =
 \sup_{\zeta\in\wt{\cX}}
 \big\{\ld\cF_R(\zeta)-\wt I(\zeta)\big\}.
\end{align}

Consider a configuration
$\zeta=\{[g_j^2(x)\,dx]\}_{j\in J}$ with finite rate.
Normalizing each $g_j$ in $L^2(\R)$ and using
\eqref{M_exp}, we obtain
\[
 2\ld\|g_j^2\|_\H
 \leq\|g_j'\|_{L^2}^2+\sfG_H(2\ld)\|g_j\|_{L^2}^2.
\]
Since $\sum_j\|g_j'\|_{L^2}^2<\infty$ and
$\sum_j\|g_j\|_{L^2}^2\leq1$, summing this inequality shows that
$\sum_j\|g_j^2\|_\H<\infty$. The Fourier representation of
$K_R$ then yields
\[
 \cF_R(\zeta)^2
 =\frac1{2\pi}\sum_j\int_{-R}^R
 |\xi|^{1-2H}|\wh{g_j^2}(\xi)|^2\,d\xi
 \leq\sum_j\|g_j^2\|_\H^2,
\]
and hence $\cF_R(\zeta)\leq\sum_j\|g_j^2\|_\H$.
Consequently,
\begin{align*}
 \ld\cF_R(\zeta)-\wt I(\zeta)
 &\leq\sum_j
 \Big(\ld\|g_j^2\|_\H-\frac12\|g_j'\|_{L^2}^2\Big) \\
 &\leq\frac12\sfG_H(2\ld)\sum_j\|g_j\|_{L^2}^2 
 \leq\frac12\sfG_H(2\ld)
 =2^{\frac{1-H}{1+H}}
 \ld^{\frac{2}{1+H}}\sfG_H.
\end{align*}
Together with \eqref{ieq_ldp}, this proves \eqref{band_lap}.

For \eqref{band_kill}, independence of $\tau$ gives
\begin{align}\label{integral_v8}
 \E\big[\exp(\ld\sqrt{\bfX_R})\big]
 =\int_0^\infty e^{-t}
 \E\big[\exp(\ld\sqrt{X_R(t)})\big]\,dt.
\end{align}
By \eqref{def_fLt-bdd},
\[
 X_R(t)\leq\frac{t^2}{2\pi}
 \int_{-R}^R|\xi|^{1-2H}\,d\xi
 \leq C_HR^{2-2H}t^2.
\]
Thus the integral in \eqref{integral_v8} is finite over
every bounded time interval.
For $0<\ld< \frac{1}{2}\ld_\ast $, \eqref{band_lap} implies
\[
 \limsup_{t\to\infty}\frac1t
 \log\E\big[\exp(\ld\sqrt{X_R(t)})\big]
 \leq(2\ld/\ld_\ast)^{2/(1+H)}<1.
\]
Thus, for some $\eps>0$ and all sufficiently large $t$,
we have 
$e^{-t}\E\big[\exp(\ld\sqrt{X_R(t)})\big] \leq e^{-\eps t}.$
The integral \eqref{integral_v8} is therefore finite, proving \eqref{band_kill}.
\end{proof}

We conclude this section with two-sided moment estimates for
$\bfY$.
For every integer $n\geq1$, the calculation in
\eqref{ieq_Y-YR-n} gives
\begin{align*}
 \E[\bfY_R^n]
 &=(2\pi)^{-n}\int_{[-R,R]^n}
 \prod_{j=1}^n|\xi_j|^{1-2H}
 \bigg[
 \sum_{\s\in\fS_n}\prod_{k=1}^n
 \frac{2}{2+\big|\sum_{\l=1}^k\xi_{\s(\l)}\big|^2}
 \bigg]^2\,d\pmxi{n}.
\end{align*}
By Lemma~\ref{lem_tail_Y}, $\bfY_R\to\bfY$ in $L^n(\O)$.
Monotone convergence on the Fourier side and
$\xi_j=\sqrt2\,\eta_j$ yield
\begin{align}\label{eq_Y-moment}
 \E[\bfY^n]
 =\lim_{R\to\infty}\E[\bfY_R^n]
 =\Big(\frac{2^{1-H}}{2\pi}\Big)^nF_n,
\end{align}
where
\begin{align}\label{def_Dn}
 F_n
 \coloneqq
 \int_{\R^n}
 \prod_{j=1}^n|\xi_j|^{1-2H}
 \bigg[
 \sum_{\s\in\fS_n}
 \prod_{k=1}^n
 \bigg(1+\Big|\sum_{\l=1}^k\xi_{\s(\l)}\Big|^2\bigg)^{-1}
 \bigg]^2\,d\pmxi{n}.
\end{align}

\begin{lemma}\label{rmk_Y-moment}
Let $-1/4<H<1/2$. There exist constants
$C_H\geq c_H>0$, depending only on $H$, such that
for every positive integer $n$,
\begin{align*} 
 c_H^n(n!)^2
 \leq\E[\bfY^n]
 \leq C_H^n(n!)^2.
\end{align*}
\end{lemma}

\begin{proof}
For the upper bound, fix $\ld=\ld_\ast/2$.
Proposition~\ref{prop_exp} and
$(2n)!\leq4^n(n!)^2$ give
\[
 \E[\bfY^n]
 \leq\E[|\bfY|^n]
 \leq\frac{(2n)!}{\ld^{2n}}
 \E\big[\exp(\ld\sqrt{|\bfY|})\big]
 \leq C_H^n(n!)^2.
\]

For the lower bound, let $Q=[-4,4]^n$ and let $h$ be the
nonnegative square root of the integrand in \eqref{def_Dn}.
Restricting the integral to $Q$ and applying Cauchy--Schwarz
to $h$ and $1$ on $Q$, whose volume is $8^n$, gives
\[
 \sqrt{F_n}
 \geq \Bigg(\int_Q h^2\,d\pmxi{n}\Bigg)^{1/2}
 \geq 8^{-n/2}\int_Q h\,d\pmxi{n}.
\]
Since every term in the permutation sum is nonnegative, $h$ is
that sum multiplied by $\prod_{j=1}^n|\xi_j|^{(1-2H)/2}$.
For each $\s\in\fS_n$, relabeling the coordinates by $\s$
preserves both $Q$ and this weight. 
Thus each of the $n!$ terms has the same integral as the term for
the identity permutation, and we obtain
\begin{align}
 \sqrt{F_n}
 &\geq \frac{n!}{8^{n/2}}
 \int_{[-4,4]^n}
 \prod_{j=1}^n
 \frac{|\xi_j|^{(1-2H)/2}}
      {1+|\xi_1+\cdots+\xi_j|^2}\,d\pmxi{n}  \notag \\
 &\geq \frac{n!}{8^{n/2}}
 \int_{R_n}
 \prod_{j=1}^n
 \frac{|\eta_j-\eta_{j-1}|^{(1-2H)/2}}
      {1+\eta_j^2}\,d\pmeta{n}, \label{transa}
\end{align}
where in the second inequality we made the change of variables
$\eta_j=\xi_1+\cdots+\xi_j$, with $\eta_0=0$, and restricted
the integral to the rectangle
\[
 R_n=\big\{(\eta_1,\ldots,\eta_n)\in\R^n \colon 
 \eta_j\in[-2,-1]\text{ for odd }j,\quad
 \eta_j\in[1,2]\text{ for even }j\big\}.
\]
The change of variables has Jacobian determinant one, and the
transformed domain is
$\{|\eta_j-\eta_{j-1}|\leq4,\ 1\leq j\leq n\}$,
which contains $R_n$.
This rectangle $R_n$ has volume one. On it, $|\eta_j-\eta_{j-1}|\geq1$
and $1+\eta_j^2\leq5$ for every $j$. Since $H<1/2$, the
transformed integrand \eqref{transa} is therefore bounded below by $5^{-n}$.
Hence, we have 
$
 \sqrt{F_n}
 \geq\frac{n!}{8^{n/2}5^n}
 \geq15^{-n}n!.
$
Combining this estimate with \eqref{eq_Y-moment} proves
the lower bound.
\end{proof}

\section{Proof of Theorem~\ref{thm_ham}: optimal solvability}\label{sec_proof}

We first show that the second chaos kernel has infinite squared norm
for every $t>0$ when $H\leq-1/4$. For $-1/4<H<1/2$, we bound the
unsymmetrized kernels in Section~\ref{ss_wker} and use that bound
in Section~\ref{ss_build} to construct the solution.

\begin{proof}[Proof of Theorem~\ref{thm_ham} \bf (nonexistence)]
Fix $t>0$ and $x\in\R$. Recall that $\wt f_{2,t,x}$ is the
symmetrization of the kernel $f_{2,t,x}$ defined in \eqref{ker_w}.
By \eqref{four_ker}, its Fourier transform depends on $x$ only
through the factor $e^{ix(\xi_1+\xi_2)}$, which has modulus one.
Thus $\|\wt f_{2,t,x}\|_{\H^{\otimes2}}^2
=\|\wt f_{2,t,0}\|_{\H^{\otimes2}}^2$, and we may take $x=0$.

Integrating the first time variable in \eqref{four_ker}, we have
$\int_0^s \wh G_r(a)\,dr=\frac{1-\cos(as)}{a^2}.$
Hence $\wh f_{2,t,0}(\xi_1,\xi_2)=K_t(\xi_1,\xi_1+\xi_2)$, where
\begin{align}\label{phi}
 K_t(a,b)
 =\int_0^t \frac{\sin((t-s)b)}{b}\,
            \frac{1-\cos(as)}{a^2}\,ds.
\end{align}

\noi
Symmetrizing in $\xi_1,\xi_2$ gives
$ \wh{\wt f}_{2,t,0}(\xi_1,\xi_2)
 =\frac12\big(K_t(\xi_1,\xi_1+\xi_2)
             +K_t(\xi_2,\xi_1+\xi_2)\big).
$
Using \eqref{cov} and changing variables to $a=\xi_1$ and
$b=\xi_1+\xi_2$, we obtain
\begin{align*}
 \|\wt f_{2,t,0}\|_{\H^{\otimes2}}^2
 &=\frac{1}{4(2\pi)^2}\int_{\R^2}
 |K_t(a,b)+K_t(b-a,b)|^2
 |a|^{1-2H}|b-a|^{1-2H}\,da\,db \\
 &\geq \frac{1}{4(2\pi)^2}\int_{D_t}
 |K_t(a,b)|^2 |a|^{1-2H}|b-a|^{1-2H}\,da\,db,
\end{align*}
where $D_t=\{(a,b)\colon a\geq4/t,\ 1/(2t)\leq b\leq1/t\}$.
The inequality follows because $0\leq(t-s)b\leq1$ on $D_t$
for every $s\in[0,t]$, so both $K_t(a,b)$ and $K_t(b-a,b)$
are nonnegative by \eqref{phi}.

For $(a,b)\in D_t$ and $0\leq s\leq t/2$, the inequality
$\sin z\geq2z/\pi$ on $[0,\pi/2]$ gives
\[
 \frac{\sin((t-s)b)}{b}
 \geq\frac{2(t-s)}{\pi}\geq\frac{t}{\pi}.
\]
Restricting the integral in \eqref{phi} to $[0,t/2]$ therefore gives
\[
 K_t(a,b)
 \geq\frac{t}{\pi a^2}\int_0^{t/2}(1-\cos(as))\,ds
 =\frac{t}{\pi a^2}
   \Big(\frac{t}{2}-\frac{\sin(at/2)}{a}\Big)
 \geq\frac{t^2}{4\pi a^2},
\]
where the last inequality uses $1/a\leq t/4$.

Finally, $3a/4\leq a-b\leq a$ on $D_t$. Since $1-2H>0$,
the weight $a^{1-2H}(a-b)^{1-2H}$ is at least
$(3/4)^{1-2H}a^{2-4H}$. Combining this with the lower bound
for $K_t(a,b)$ and integrating in $b$ over $[1/(2t),1/t]$,
we obtain, for some $C_{t,H}>0$,
\begin{align*}
 \|\wt f_{2,t,0}\|_{\H^{\otimes2}}^2
 &\geq C_{t,H}\int_{4/t}^{\infty}a^{-4}a^{2-4H}\,da
 = C_{t,H}\int_{4/t}^{\infty}a^{-2-4H}\,da
 =\infty
 \quad\text{if } H\leq-1/4.
\end{align*}
Since this holds for every $t>0$, the necessity part of
Proposition~\ref{prop_chaos} rules out a square-integrable Skorokhod random-field solution
on any interval $[0,T)$ with $T>0$.
\end{proof}

\subsection{Chaos coefficients and the convolution integrand}
\label{ss_wker}

The chaos expansion of the solution uses $f_{n,t,x}$ symmetrized
in all $n$ spatial variables. In the time convolution, only the
first $n-1$ variables are symmetrized, as the last variable
is reserved in the Skorokhod integral. Both symmetrizations are
contractions in $\H^{\otimes n}$, so a bound on $f_{n,t,x}$
controls both kernels.

\begin{lemma}\label{lem_wker}
Let $-1/4<H<1/2$. There exists $C_H>0$ such that, for every
$n\geq1$ and $(t,x)\in\R_+\times\R$,
\begin{equation}\label{unsym_resolvent_bd}
 \|f_{n,t,x}\|_{\H^{\otimes n}}^2
 \leq \frac{C_H^n t^{(2H+2)n}}{\Gam((2H+2)n+1)}.
\end{equation}
In particular,
\begin{equation}\label{wker_bd}
 \|\wt f_{n,t,x}\|_{\H^{\otimes n}}^2
 \leq \|f_{n,t,x}\|_{\H^{\otimes n}}^2
 \leq \frac{C_H^n t^{(2H+2)n}}{(n!)^{2H+2}}.
\end{equation}
\end{lemma}

Balan and Yuan \cite{BY24} proved \eqref{wker_bd} for
$1/4<H<1/2$. They first fix the time increments and integrate
the squared Fourier transform of each wave kernel in the spatial
frequency. For a time increment $r>0$, the resulting integrals
have the form
$\int_\R |\eta|^\al\,\frac{\sin^2(r\eta)}{\eta^2}\,d\eta.$
These integrals are finite at infinity precisely when $\al<1$.
Since the largest power contributed by the spatial weights is
$\al=2(1-2H)$, this step requires $H>1/4$.
We extend the bound to $-1/4<H<1/2$ by first taking
the Laplace transform in the terminal time $t$. The resulting
wave resolvent also appears in \cite[Proposition~3.5]{CDST25}.
Its squared modulus is integrable against $|\eta|^\al\,d\eta$ for
$0\leq\al<3$, allowing the largest spatial exponent $2(1-2H)$
precisely when $H>-1/4$.

\begin{proof}[Proof of Lemma~\ref{lem_wker}]
Set $a=1-2H\in(0,3/2)$.
We first show that, for $0\leq\al<3$ and $\o\in\R$,
\begin{equation}\label{scalar_resolvent_bd}
 \int_\R \frac{|\eta|^\al}{|(1-i\o)^2+\eta^2|^2}\,d\eta
 \leq C_\al(1+\o^2)^{(\al-2)/2}.
\end{equation}

The factorization
\[
 |(1-i\o)^2+\eta^2|^2
 =|1+i(\eta-\o)|^2\,|1-i(\eta+\o)|^2
 =[1+(\eta-\o)^2][1+(\eta+\o)^2]
\]
shows that the integrand is even in both $\eta$ and $\o$.
We may therefore take $\o\geq0$ and double the integral over $\eta\geq0$.
Set $u=\eta-\o$. Since $\eta\leq\o+|u|$, we have
$1\leq(1+\eta+\o)/(1+\o)\leq2+|u|$. Hence
\[
 \frac{\eta^\al}{1+(\eta+\o)^2}
 \leq 2(1+\eta+\o)^{\al-2}
 \leq 2(1+\o)^{\al-2}(2+|u|)^{(\al-2)_+}.
\]
The last inequality uses the lower bound on the ratio when $\al\leq2$
and the upper bound when $\al>2$. Since $(\al-2)_+<1$, integrating
against $du/(1+u^2)$ gives \eqref{scalar_resolvent_bd}.

The spatial phase in \eqref{four_ker} has modulus one, so we take $x=0$
and set $\eta_j=\xi_1+\cdots+\xi_j$, with $\eta_0=0$.
To take the Laplace transform in $t$, set $s_{n+1}=t$, $r_0=s_1$,
and $r_j=s_{j+1}-s_j$ for $1\leq j\leq n$. Then
$t=r_0+\cdots+r_n$, and the change of variables maps the region
$0<s_1<\cdots<s_n<t$ onto $(0,\infty)^{n+1}$ with Jacobian one.
Thus, for $z=1-i\o$, \eqref{four_ker} gives
\begin{align*}
 \int_0^\infty e^{-zt}g_{n,t,0}(\pmxi{n})\,dt
 &=\int_0^\infty e^{-zt}
   \int_{0<s_1<\cdots<s_n<t}
   \prod_{j=1}^n\wh G_{s_{j+1}-s_j}(\eta_j)
   \,ds_1\cdots ds_n\,dt \\
 &=\int_{(0,\infty)^{n+1}}
   e^{-z(r_0+\cdots+r_n)}
   \prod_{j=1}^n\wh G_{r_j}(\eta_j)\,dr_0\cdots dr_n \\
 &=\bigg(\int_0^\infty e^{-zr_0}\,dr_0\bigg)
   \prod_{j=1}^n
   \bigg(\int_0^\infty e^{-zr_j}\wh G_{r_j}(\eta_j)\,dr_j\bigg)
 =\frac1z\prod_{j=1}^n\frac1{z^2+\eta_j^2},
\end{align*}
where the last equality follows from \eqref{Lapa}.

Since $|\wh G_r(\eta)|\leq r$, we have
\[
 |g_{n,t,0}(\pmxi{n})|
 \leq\int_{0<s_1<\cdots<s_n<t}
 \prod_{j=1}^n(s_{j+1}-s_j)\,ds_1\cdots ds_n
 =\frac{t^{2n}}{(2n)!}.
\]
Hence, for each fixed frequency vector, $e^{-t}g_{n,t,0}(\pmxi{n})$,
extended by zero to $t<0$, belongs to $L^2(\R)$. Its Fourier
transform in $t$ is the Laplace transform computed above,
so Plancherel's identity gives
\[
 \int_0^\infty e^{-2t}|g_{n,t,0}(\pmxi{n})|^2\,dt
 =\frac1{2\pi}\int_\R\frac{d\o}{1+\o^2}
   \prod_{j=1}^n\frac1{|(1-i\o)^2+\eta_j^2|^2}.
\]
We now integrate against
$(2\pi)^{-n}\prod_{j=1}^n|\xi_j|^a\,d\pmxi{n}$ and apply
Tonelli's theorem. Since the change of variables
$\xi_j=\eta_j-\eta_{j-1}$ has Jacobian one, we obtain
\begin{align}\label{inner1}
 \int_0^\infty e^{-2t}\|f_{n,t,0}\|_{\H^{\otimes n}}^2\,dt
 &=\frac1{(2\pi)^{n+1}}
 \int_\R\frac{d\o}{1+\o^2}
 \int_{\R^n}\prod_{j=1}^n
 \frac{|\eta_j-\eta_{j-1}|^a}
      {|(1-i\o)^2+\eta_j^2|^2}\,d\pmeta{n}.
\end{align}

To bound the inner integral, we use the same expansion as in
\eqref{ieq_Y-YR}. Applying
\begin{align*}
 |\eta_j-\eta_{j-1}|^a
 \leq C_H\big(|\eta_j|^a+|\eta_{j-1}|^a\big)
\end{align*}
to each factor and expanding the product bounds the numerator by
$C_H^n$ times a sum of at most $2^n$ monomials. Each nonzero
monomial has the form $\prod_{j=1}^n|\eta_j|^{\al_j}$ with
$\pmal{n}\in\fA_n$, where $\fA_n$ is defined in \eqref{def_A-n}.
In particular, $\sum_{j=1}^n\al_j=an$ and $\al_j\leq 2a<3$
for every $j$, since $H>-1/4$.

For each such monomial, \eqref{scalar_resolvent_bd} gives

\noi
\begin{align*}
 \prod_{j=1}^n
 \int_\R\frac{|\eta_j|^{\al_j}}
 {|(1-i\o)^2+\eta_j^2|^2}\,d\eta_j
 \leq C_H^n(1+\o^2)^{\frac12\sum_{j=1}^n(\al_j-2)}
 =C_H^n(1+\o^2)^{-(1+2H)n/2}. 
\end{align*}

\noi
Summing over the monomials and absorbing the factor $2^n$ into
$C_H^n$ gives the same bound for the inner integral in \eqref{inner1}.
Since $1+2H>0$, integration against $(1+\o^2)^{-1}\,d\o$ gives
\[
  \int_0^\infty e^{-2t}\|f_{n,t,0}\|_{\H^{\otimes n}}^2\,dt
 \leq C_H^n.
\]

For $t>0$, rescaling the time variables in \eqref{four_ker} gives
\begin{equation}\label{ker_scale}
 g_{n,t,x}(\pmxi{n})
 =e^{ix(\xi_1+\cdots+\xi_n)}t^{2n}g_{n,1,0}(t\pmxi{n}).
\end{equation}
Changing variables in the spatial Fourier norm yields
$\|f_{n,t,x}\|_{\H^{\otimes n}}^2
=t^{(2H+2)n}\|f_{n,1,0}\|_{\H^{\otimes n}}^2$.
Substitution into the preceding time integral gives
\[
 \frac{\Gam((2H+2)n+1)}{2^{(2H+2)n+1}}\,
 \|f_{n,1,0}\|_{\H^{\otimes n}}^2\leq C_H^n.
\]
Together with scaling, this proves \eqref{unsym_resolvent_bd} for
$t>0$. At $t=0$ the kernel vanishes. Finally, \eqref{gam_low} gives
$\Gam((2H+2)n+1)\geq c_{2H+2}^{\,n}(n!)^{2H+2}$, and symmetrization
is a contraction by \eqref{sym_bd}. This proves \eqref{wker_bd}.
\end{proof}

\subsection{Construction of the solution}\label{ss_build}

We use the standard Wiener chaos construction and include the proof
for completeness. The kernel estimate \eqref{wker_bd} gives both
convergence of the chaos series and the moment bound \eqref{ham_up}.
Two points require care when $-1/4<H\leq0$. First, the zero initial
velocity holds after testing against Schwartz functions, while the
first chaos prevents the pointwise difference quotient from converging
to zero in $L^2(\O)$ (see Remark~\ref{rem_vel}). Second, $G_r\notin\H$
for every $r>0$. We therefore integrate in time first and use
\eqref{wker_bd} for the partially symmetrized kernels to verify that
the resulting convolution belongs to the domain of the  
Skorokhod integral.

\begin{proposition}\label{prop_ham_up}
Let $-1/4<H<1/2$.
Then \eqref{cond_w} holds for every $(t,x)\in\R_+\times\R$,
and equation \eqref{ham} has a unique global random-field solution
$u$. Moreover, there exists $C_H\in(0,\infty)$, depending only
on $H$, such that, for every $p\in [2, \infty)$ and
$(t,x)\in\R_+\times\R$,
\begin{align}\label{ham_up}
 \|u(t,x)\|_p
 \leq C_H\exp\!\big\{C_H
 (p-1)^{\frac{1}{2H+1}}t^{\frac{2H+2}{2H+1}}\big\}.
\end{align}
\end{proposition}

\begin{proof}
Set $u^{[N]}(t,x)=\sum_{n=0}^N I_n(\wt f_{n,t,x})$. For
$2\leq p<\infty$ and $M>N$, hypercontractivity \eqref{chaos_hyp}
and the kernel bound \eqref{wker_bd} give
\begin{align}\label{ham_tail}
 \|u^{[M]}(t,x)-u^{[N]}(t,x)\|_p
 &\leq\sum_{n=N+1}^M (p-1)^{n/2}\sqrt{n!}\,
       \|\wt f_{n,t,x}\|_{\H^{\otimes n}}\notag\\
 &\leq\sum_{n=N+1}^M
 \frac{[C_H\sqrt{p-1}\,t^{H+1}]^n}{(n!)^{H+1/2}}.
\end{align}
Since $H+\frac{1}{2}>0$, the sums $u^{[N]}$ converge in $L^p(\O)$ to a
common limit $u$, uniformly on $[0,T]\times\R$ for every finite $T$.
Orthogonality at $p=2$ gives \eqref{cond_w}. Taking $N=0$ and
$M\to\infty$ in \eqref{ham_tail} and using \eqref{series} yields
\[
 \|u(t,x)\|_p
 \leq 1+\sum_{n=1}^\infty
 \frac{[C_H\sqrt{p-1}\,t^{H+1}]^n}{(n!)^{H+1/2}}
 \leq C_H\exp\!\big\{C_H
 (p-1)^{\frac1{2H+1}}t^{\frac{2H+2}{2H+1}}\big\}.
\]
This proves \eqref{ham_up} and, in particular, bounds $\|u(t,x)\|_2$
uniformly on each such strip.

For $t>0$, \eqref{ker_scale} gives joint continuity of
$(t,x)\mapsto f_{n,t,x}$ in $\H^{\otimes n}$, by strong continuity
of dilations and phase multiplication in the weighted Fourier $L^2$
space. At $t=0$, \eqref{wker_bd} gives continuity uniformly in $x$.
By \eqref{iso} and uniform convergence of the chaos series, $u$ is
jointly $L^2(\O)$ continuous. On each $[0,m]\times[-m,m]$, choose a
piecewise constant grid approximation with uniform $L^2(\O)$ error
at most $2^{-m}$, using $\s\{W\}$-measurable values at the grid
points. Extend it by zero outside the rectangle. At every fixed
$(t,x)$ these approximations converge almost surely to $u(t,x)$.
Their limit, set to zero where it does not exist, gives a jointly
measurable, $\s\{W\}$-measurable version. Spatial stationarity
follows from the stationary noise and translation covariance of
the kernels.

We check the initial conditions in Definition~\ref{def_sol}(i).
At $t=0$ only the zeroth chaos remains, so $u(0,x)=1$.
For $0<t\leq1$, taking $p=2$ in \eqref{ham_tail} gives
\begin{equation}\label{init_pos}
 \sup_{x\in\R}\|u(t,x)-1\|_2\leq C_Ht^{H+1}.
\end{equation}
Starting the same sum at $n=2$ bounds the remaining chaoses by
$C_Ht^{2H+2}$ in $L^2(\O)$, uniformly in $x$. For
$\varphi\in\cS(\R)$ and $t>0$, the pairing
$\int_\R\varphi(x)(u(t,x)-1)\,dx$ exists as a Bochner integral
in $L^2(\O)$ by continuity and \eqref{ham_up}. By \eqref{four_ker},
the Fourier transform of the paired first-chaos kernel is
$\wh\varphi(\xi)(1-\cos(t\xi))/\xi^2$. Since
$1-\cos(t\xi)\leq t^2\xi^2/2$, \eqref{iso} bounds the $L^2(\O)$
norm of this chaos by $(t^2/2)\|\varphi\|_\H$.
Adding the higher chaoses gives
\begin{equation}\label{init_vel}
 \bigg\|\frac1t\int_\R\varphi(x)(u(t,x)-1)\,dx\bigg\|_2
 \leq\frac t2\|\varphi\|_\H+C_Ht^{2H+1}\|\varphi\|_{L^1(\R)}
\to 0
\quad\text{as $t\downarrow 0$.}
\end{equation}

For conditions~(ii) and~(iii) in Definition~\ref{def_sol}, fix $t>0$ and $x\in\R$. Set
\[
 \Psi^{[N]}_{t,x}(y)=\int_0^t G_{t-s}(x-y)u^{[N]}(s,y)\,ds
\AND
 \Psi_{t,x}(y)=\int_0^t G_{t-s}(x-y)u(s,y)\,ds.
\]
For each $y$, these are Bochner integrals in $L^2(\O)$ since the
random fields are $L^2(\O)$ continuous and the wave kernel is bounded
and measurable. By \eqref{ker_w}, integration against the last wave
kernel gives the $n$th chaos kernel of $\Psi^{[N]}_{t,x}$, $n\leq N$,
as $f_{n+1,t,x}$ symmetrized in its first $n$ variables. Denote it
by $k_{n;t,x}$. Partial symmetrization is a contraction, so
\eqref{wker_bd} gives
\begin{equation}\label{conv_sum}
 \sum_{n=0}^{\infty}n!\|k_{n;t,x}\|_{\H^{\otimes(n+1)}}^2
 \leq\sum_{n=0}^{\infty}
 \frac{n!\,[C_Ht^{2H+2}]^{n+1}}{((n+1)!)^{2H+2}}<\infty.
\end{equation}
Hence $\Psi^{[N]}_{t,x}$ has a limit $V_{t,x}$ in $L^2(\O;\H)$.
On the other hand, $\int_\R G_{t-s}(x-y)\,dy=t-s$, so
\[
 \|\Psi^{[N]}_{t,x}-\Psi_{t,x}\|_{L^2(\O;L^1(\R))}
 \leq\frac{t^2}{2}
 \sup_{0\leq s\leq t,\,y\in\R}\|u^{[N]}(s,y)-u(s,y)\|_2
 \longrightarrow0.
\]
To identify the limits, write $\mu(d\xi)=(2\pi)^{-1}|\xi|^{1-2H}\,d\xi$.
On each bounded interval $I$, the Fourier transforms of
$\Psi^{[N]}_{t,x}$ converge to $\wh V_{t,x}$ in
$L^2(\O\times I,\PP\otimes\mu)$ by \eqref{conv_sum}.
They also converge there to $\wh\Psi_{t,x}$ by the preceding
estimate, since $|\wh h(\xi)|\leq\|h\|_{L^1}$ and $\mu(I)<\infty$.
Thus $\Psi_{t,x}$ represents $V_{t,x}$ in $L^2(\O;\H)$.

Full symmetrization gives $\wt k_{n;t,x}=\wt f_{n+1,t,x}$.
By \eqref{div}, we have 
$\dl(\Psi^{[N]}_{t,x})=u^{[N+1]}(t,x)-1$,
which converges to  $u(t,x)-1$ in $L^2(\O)$.
Closedness of $\dl$ now gives $\Psi_{t,x}\in\Dom\dl$ and
$\dl(\Psi_{t,x})=u(t,x)-1$. This verifies conditions~(ii) and~(iii)
and completes the construction. Uniqueness follows from the kernel
identification in the proof of Proposition~\ref{prop_chaos} below.
\end{proof}

\begin{proof}[Proof of Proposition~\ref{prop_chaos}]
Fix $H<1/2$ and $T\in(0,\infty]$.

\emph{Sufficiency. Suppose \eqref{cond_w} holds on $[0,T)\times\R$.
Fix $0<t<T$ and take $x=0$. The $n=2$ term in \eqref{cond_w} is finite.
The lower bound derived from \eqref{phi} then gives
$\int_{4/t}^\infty a^{-2-4H}\,da<\infty$, so $H>-1/4$.
The global solution constructed in the proof of Proposition~\ref{prop_ham_up}
restricts to a solution on $[0,T)\times\R$.}

\emph{Necessity.} Let $u$ satisfy Definition~\ref{def_sol}, with
$n$-th chaos kernel $a_{n;t,x}\in\H^{\odot n}$. For a real Schwartz function
$V$ satisfying $\wh V\in C_c^\infty(\R\setminus\{0\})$, set
$h=(-\Dl)^{H-1/2}V$ and $U_n(t,x)=\E[u(t,x)I_n(h^{\otimes n})]/n!$.
Then $h\in\H$ by \eqref{E_H'} and $U_0=1$ by \eqref{mild_w}.
For $0<t<T$, Definition~\ref{def_sol} (ii) gives $\Psi_{t,x}\in\Dom\dl$.
By the divergence formula \eqref{div} \cite[Proposition~1.3.7]{Nua06},
symmetrizing the chaos kernel of order $n-1$ of $\Psi_{t,x}$ gives $a_{n;t,x}$.
Pairing with $h^{\otimes n}$ and using \eqref{iso} and the identity
$\langle\Psi_{t,x},h\rangle_\H=\int_\R\Psi_{t,x}(y)V(y)\,dy$ from
\eqref{id_lap-phi-psi}, we obtain
\[
 \begin{aligned}
 U_n(t,x)
 &=\frac{1}{(n-1)!}\E\!\left[
     \langle\Psi_{t,x},h\rangle_\H\,I_{n-1}(h^{\otimes(n-1)})\right]\\
 &=\int_0^t\!\int_\R G_{t-s}(x-y)V(y)
     \frac{\E[u(s,y)I_{n-1}(h^{\otimes(n-1)})]}{(n-1)!}\,dy\,ds\\
 &=\int_0^t\!\int_\R G_{t-s}(x-y)V(y)U_{n-1}(s,y)\,dy\,ds,
 \qquad n\geq1.
 \end{aligned}
\]

By \eqref{iso},
$U_n(t,x)=\langle a_{n;t,x},h^{\otimes n}\rangle_{\H^{\otimes n}}$.
With $(s_{n+1},y_{n+1})=(t,x)$, iterating the recursion from $U_0=1$ gives
\[
 \begin{aligned}
 U_n(t,x)
 &=\int_{\R^n}\left[\int_{0<s_1<\cdots<s_n<t}
   \prod_{j=1}^n G_{s_{j+1}-s_j}(y_{j+1}-y_j)\,d\pms{n}\right]
   \prod_{j=1}^n V(y_j)\,d\pmy{n}\\
 &=\int_{\R^n}f_{n,t,x}(\pmy{n})\prod_{j=1}^n V(y_j)\,d\pmy{n}
 =\int_{\R^n}\wt f_{n,t,x}(\pmy{n})\prod_{j=1}^n V(y_j)\,d\pmy{n}.
 \end{aligned}
\]
The inner integral is $f_{n,t,x}(\pmy{n})$ by \eqref{ker_w}. The last
equality follows by averaging over permutations of $y_1,\ldots,y_n$,
since $\prod_{j=1}^n V(y_j)$ is symmetric.
By polarization, complex linearity and \eqref{id_lap-phi-psi},
$\wh a_{n;t,x}$ and $\wh{\wt f}_{n,t,x}$ have the same integrals against
all products of test functions in $C_c^\infty(\R\setminus\{0\})$ and
hence agree almost everywhere on $(\R\setminus\{0\})^n$. This set has
full product spectral measure, so $\wt f_{n,t,x}$ represents
$a_{n;t,x}$ in $\H^{\odot n}$.

By \eqref{iso},
\[
 \sum_{n=1}^\infty n!\|\wt f_{n,t,x}\|_{\H^{\otimes n}}^2
 =\E[|u(t,x)|^2]-1<\infty .
\]
At $t=0$ the kernels in \eqref{ker_w} vanish. This proves \eqref{cond_w}.

\emph{Uniqueness.} The kernel identification gives the same chaos
expansion \eqref{chaos-u} for every solution, with convergence in
$L^2(\O)$ by \eqref{iso}. Hence any two solutions agree in $L^2(\O)$
at each fixed $(t,x)\in[0,T)\times\R$.
\end{proof}

\begin{remark}[Pointwise initial velocity]
\rm
\label{rem_init}\label{rem_vel}
For $H>0$, \eqref{init_pos} gives
\[
 \sup_{x\in\R}\Big\|\frac{u(t,x)-1}{t}\Big\|_2
 \leq C_Ht^H\longrightarrow0
 \quad\text{as $t\downarrow0$.}
\]
For $-1/4<H\leq0$, orthogonality of the chaos and the first-chaos
formula give
\[
 \Big\|\frac{u(t,x)-1}{t}\Big\|_2^2
 \geq \frac{t^{2H}}{2\pi}
 \int_\R\frac{(1-\cos\xi)^2}{|\xi|^{3+2H}}\,d\xi.
\]
The integral is strictly positive and finite. Its integrand is
$O(|\xi|^{1-2H})$ near zero and $O(|\xi|^{-3-2H})$ at infinity.
The lower bound is a positive constant for $H=0$ and tends to
infinity for $H<0$. Thus, in this range, the spatially tested
limit \eqref{init_vel} holds even though the pointwise
$L^2(\O)$ limit to zero fails.
\end{remark}

\begin{proof}[Proof of Theorem~\ref{thm_ham} \bf(existence)]
For $-1/4<H<1/2$, Proposition~\ref{prop_ham_up} gives a unique
global solution in the sense of Definition~\ref{def_sol}.
\end{proof}

\section{Proof of Theorem~\ref{thm_mom}: sharp lower and upper moment bounds}\label{sec_mom}

Fix $-1/4<H<1/2$. The moment bound \eqref{ham_up} in
Proposition~\ref{prop_ham_up} gives the scale
\[
 t_{p,H}=(p-1)^{1/(2H+1)}t^{(2H+2)/(2H+1)}.
\]
The next proposition identifies the leading coefficient on this
scale and proves the limit in Theorem~\ref{thm_mom}. Recall from
\eqref{tpHK} that
\begin{align}\label{KH2}
 \sfK_H=\frac{2H+1}{2}
 \bigg(\frac{\sqrt{\sfG_H}}{1+H}\bigg)^{(2H+2)/(2H+1)}.
\end{align}
The   constant $\sfG_H$ in \eqref{G_H} enters both the
optimization over deterministic potentials for the lower bound
and the Brownian exponential moment estimate for the upper bound.

\begin{proposition}\label{prop_mom}
For every $\eps>0$, there exists $C_{H,\eps}\geq1$ such that,
for all $(t,x)\in\R_+\times\R$ and $p\in[2,\infty)$, the solution
$u$ to \eqref{ham} satisfies
\begin{equation}\label{mom_bounds}
 \frac{1}{2}\exp\!\big(\sfK_H t_{p,H}\big)
 \leq \|u(t,x)\|_p
 \leq C_{H,\eps}\exp\!\big((\sfK_H+\eps)t_{p,H}\big).
\end{equation}
Consequently, for every $x\in\R$,  we have 
\begin{align}\label{joint_lim}
\text{$t_{p,H}^{-1} \log\|u(t,x)\|_p  \to  \sfK_H$
as $t_{p,H}\to\infty$.}
\end{align}
\end{proposition}

Section~\ref{ss_lower} derives the lower bound from the growth of a
deterministic wave with potential $V$. Section~\ref{ss_upper} uses
the Brownian moment formula \eqref{eq_Y-moment} to estimate the
symmetrized chaos kernels at time one and sums their contributions
to obtain the upper bound. Section~\ref{sec_aux} proves the auxiliary
estimates used in the upper-bound argument.

\subsection{Sharp lower bound}\label{ss_lower}

To prove the lower bound, we test $u(t,x)$ against a Wick exponential,
as in \cite[Proposition~2.1]{Che17} for the parabolic Anderson model
and \cite[Proposition~6.3]{BCC22} for the wave equation with
time-independent noise. We take potentials $V\in\mathring\H^*$
as defined in \eqref{Hast2}. Schwartz regularity allows us to apply
the deterministic wave estimates in Section~\ref{ss_test}, and the
 norm \eqref{Hast} gives the cost of the test.
For $H>0$, every real Schwartz potential has finite norm. When
$H\leq0$, the weight $|\xi|^{2H-1}$ in \eqref{Hast} and continuity
of $\wh V$ force $\wh V(0)=\int_\R V(y)\,dy=0$. Conversely, if
$\wh V(0)=0$, then $\wh V(\xi)=O(|\xi|)$ near zero, and the norm is
finite since $2H+1>-1$. Thus, when $H\leq0$, every nonzero
$V\in\mathring\H^*$ takes both positive and negative values.

Recall the variational definition \eqref{spec_val} of $\Ld(V)$
and its spectral characterization \eqref{spec_bot}. When
$\Ld(V)>0$, the positive ground state gives the lower bound
\eqref{wave_sup} for the deterministic wave $w_V$. We transfer
this bound to the moments of $u$.

\begin{lemma}\label{lem_prof}
Let $V\in\mathring\H^*$ satisfy $\Ld(V)>0$. Then for every $t>0$,
$x\in\R$, and $p\in[2,\infty)$,
\begin{equation}\label{prof_bd}
 \log\|u(t,x)\|_p
 \geq t_{p,H}\Big(\sqrt{\Ld(V)}-\frac{\|V\|_{\H^*}^2}{2}\Big)
       -\log 2.
\end{equation}
\end{lemma}

\begin{proof}
Set $h_V=(-\Dl)^\al V$, where $\al=H-\frac12$. By
\eqref{id_lap-phi-psi}, $h_V$ represents pairing against $V$ in
$\H$, and \eqref{E_H'} gives $\|h_V\|_\H=\|V\|_{\H^*}$. The Wick
exponential $\cZ=\exp(W(h_V)-\tfrac12\|h_V\|_\H^{2})$
has the chaos expansion
\begin{equation}\label{chaos-z}
 \cZ=\sum_{n=0}^\infty\frac{1}{n!}I_n(h_V^{\otimes n})
 \quad\text{in }L^2(\O)
\end{equation}
by \cite[Proposition~3.5 and Example~5.6]{Hu17}. For $p'=p/(p-1)$,
the Gaussian moment formula gives
\begin{equation}\label{exp_norm}
 \|\cZ\|_{p'}
 =\exp\!\Big(\frac{p'-1}{2}\|h_V\|_\H^2\Big)
 =\exp\!\Big(\frac{\|V\|_{\H^*}^2}{2(p-1)}\Big).
\end{equation}

Since both chaos series \eqref{chaos-u} and \eqref{chaos-z} converge
in $L^2(\O)$, orthogonality \eqref{iso} and the symmetry of
$h_V^{\otimes n}$ give the first equality below. The second follows
by identifying each summand with the corresponding Picard term.
\begin{equation}\label{S_trans}
 \E[u(t,x)\cZ]
 =\sum_{n=0}^\infty
   \langle f_{n,t,x},h_V^{\otimes n}\rangle_{\H^{\otimes n}}
 =w_V(t,x),
\end{equation}
where the term at $n=0$ is $1$. For $n\geq1$, the kernel
$f_{n,t,x}$ is bounded and supported in $[x-t,x+t]^n$ by
\eqref{ker_w}. Applying duality and Plancherel in each variable gives
\[
 \langle f_{n,t,x},h_V^{\otimes n}\rangle_{\H^{\otimes n}}
 =\int_{\R^n}f_{n,t,x}(\bsy_n)\prod_{j=1}^n V(y_j)\,\rmd\bsy_n,
\]
which is the term of order $n$ in the Picard series \eqref{wave_ser}.

By Lemma~\ref{wave_grow}, H\"older's inequality, and spatial
stationarity, we deduce from \eqref{S_trans} and \eqref{exp_norm} that
\begin{equation}\label{holder_sup}
 \cosh\big(t\sqrt{\Ld(V)}\big)
 \leq \sup_{y\in\R}|w_V(t,y)|
 \leq \|u(t,x)\|_p
       \exp\!\Big(\frac{\|V\|_{\H^*}^2}{2(p-1)}\Big).
\end{equation}

Now for $\ld>0$, we set $V_\ld(y)=\ld^2V(\ld y)$.
Changing variables $\xi=\ld\eta$ in \eqref{Hast} and using the
$L^2$-preserving dilation $g_\ld(y)=\ld^{1/2}g(\ld y)$ in
\eqref{spec_val}, we obtain
\[
 \|V_\ld\|_{\H^*}^2=\ld^{2H+2}\|V\|_{\H^*}^2
 \AND  \Ld(V_\ld)=\ld^2\Ld(V).
\]
Thus $V_\ld\in\mathring\H^*$ and $\Ld(V_\ld)>0$. Applying
\eqref{holder_sup} to $V_\ld$ and using $\cosh r\geq e^r/2$ for
$r\geq0$, we obtain
\[
 \log\|u(t,x)\|_p
 \geq \ld t\sqrt{\Ld(V)}
 -\frac{\ld^{2H+2}}{2(p-1)}\|V\|_{\H^*}^2-\log 2.
\]
Taking $\ld=[(p-1)t]^{1/(2H+1)}$, we have
$\ld^{2H+2}=(p-1)\ld t$ and $\ld t=t_{p,H}$.
Substituting these identities into the preceding bound yields
\eqref{prof_bd}.
\end{proof}

We next optimize the right-hand side of \eqref{prof_bd} over
all $V\in\mathring\H^*$ with $\Ld(V)>0$.

\begin{lemma}\label{common_C}
For $-1/4<H<1/2$, with $\sfK_H$ as in \eqref{KH2},
we have
\[
 \sup
 \Big\{\sqrt{\Ld(V)}-\frac{\|V\|_{\H^*}^2}{2} \,\colon \, V\in\mathring\H^*  \,\, {\rm and} \,\,
  \Ld(V)>0\Big\}=\sfK_H.
\]
\end{lemma}

\begin{proof}
Fix $\kk>0$ and restrict to potentials with $\|V\|_{\H^*}^2=\kk$.
For $g\in W^{1,2}(\R)$, Lemma~\ref{lem_prod} gives $g^2\in\H$.
By covariance duality \eqref{id_lap-phi-psi} and density of
$\mathring\H^*$ on the sphere of radius $\sqrt{\kk}$ in $\H^*$,
we have
\begin{equation}\label{dual_id}
 \sup_{\substack{V\in\mathring\H^*\,,\ \|V\|_{\H^*}^2=\kk}}
 \int_\R V(x)g(x)^2\,dx=\sqrt{\kk}\,\|g^2\|_\H.
\end{equation}
Interchanging the suprema over $V$ and $g$ in \eqref{spec_val}
and using \eqref{dual_id} and \eqref{M_exp}, we obtain
\[
 \sup_{\substack{V\in\mathring\H^*\\ \|V\|_{\H^*}^2=\kk}}\Ld(V)
 =\sup_{\substack{g\in W^{1,2}(\R)\\ \|g\|_{L^2}=1}}
   \big\{\sqrt{\kk}\,\|g^2\|_\H-\|g'\|_{L^2}^2\big\}
 =\sfG_H(\sqrt{\kk})=\kk^{1/(1+H)}\sfG_H>0.
\]
Since this supremum is positive, it can be approached by potentials
with $\Ld(V)>0$. Thus the supremum in Lemma~\ref{common_C},
restricted to $\|V\|_{\H^*}^2=\kk$, equals
$\sqrt{\sfG_H}\,\kk^{1/(2+2H)}-\frac{1}{2}\kk.$
Optimizing over $\kk>0$, we find the unique maximizer
\[
 \kk_*=\bigg(\frac{\sqrt{\sfG_H}}{1+H}\bigg)^{(2+2H)/(1+2H)}.
\]
At $\kk_*$, the first term equals $(1+H)\kk_*$, so the maximum is
$\frac{1+2H}{2}\kk_*=\sfK_H$ by \eqref{tpHK}.
\end{proof}

\begin{proof}[\bf \raisebox{0.15ex}{\(\bullet\)}  Proof of the lower bound in Proposition~\ref{prop_mom}]
For $t>0$, take the supremum over $V$ in \eqref{prof_bd} and use
Lemma~\ref{common_C}. At $t=0$, the lower bound in \eqref{mom_bounds}
follows from $u(0,x)=1$.
\end{proof}

\subsection{Sharp upper bound}\label{ss_upper}

The upper bound in Proposition~\ref{prop_mom} follows by summing
the chaos series using the estimate \eqref{ieq_claim} below.
This refines \eqref{wker_bd} from Section~\ref{ss_wker} by
giving the exponential rate needed for the sharp moment bound.

\begin{lemma}\label{claim}
Let $-1/4<H<1/2$. For every $\dl>0$ there exists $C_{H,\dl}>0$
such that, for all $n\geq0$,
\begin{equation}\label{ieq_claim}
 \|\wt f_{n,1,0}\|_{\H^{\otimes n}}^2
 \leq \frac{C_{H,\dl}}{(n!)^{2H+2}}
 \Big(\frac{2\sfK_H}{1+2H}+\dl\Big)^{(1+2H)n}.
\end{equation}
\end{lemma}

\begin{proof}[\bf \raisebox{0.15ex}{\(\bullet\)} Proof of the upper bound in Proposition~\ref{prop_mom}]
By spatial stationarity, it suffices to take $x=0$. Fix $t>0$.
Since dilation commutes with symmetrization, the change of variables
following \eqref{ker_scale} also gives
\begin{equation}\label{coeff_scale}
 \|\wt f_{n,t,0}\|_{\H^{\otimes n}}^2
 =t^{(2H+2)n}\|\wt f_{n,1,0}\|_{\H^{\otimes n}}^2,
 \quad t>0,\quad n\geq0.
\end{equation}
Set $a=H+\frac12$ and recall $t_{p,H}=(p-1)^{1/(2H+1)}t^{(2H+2)/(2H+1)}$. Thus
$t_{p,H}^a=\sqrt{p-1}\,t^{H+1}$. Given $\eps>0$, choose $\dl>0$
so that $a\dl<\eps$. Applying hypercontractivity \eqref{chaos_hyp}, and using \eqref{coeff_scale} and
\eqref{ieq_claim}, we obtain
\[
 \|u(t,0)\|_p
 \leq \sum_{n=0}^\infty \sqrt{n!}\,t_{p,H}^{an}
      \|\wt f_{n,1,0}\|_{\H^{\otimes n}}
 \leq \sqrt{C_{H,\dl}}\sum_{n=0}^\infty
       \frac{[( \frac{\sfK_H}{a}+\dl)t_{p,H}]^{an}}{(n!)^a}.
\]
By \eqref{series-asymp}, the logarithm of the last series is
$(\sfK_H+a\dl)t_{p,H}+o(t_{p,H})$ as $t_{p,H}\to\infty$. 
Since $a\dl<\eps$ and the series is bounded on compact intervals,
we obtain
$\|u(t,0)\|_p\leq C_{H,\eps}e^{(\sfK_H+\eps)t_{p,H}}$
for all $(t,p)\in (0,\infty)\times [2, \infty)$.
The bound extends to $t=0$ since $u(0,0)=1$, and we may take
$C_{H,\eps}\geq1$. Together with the lower bound in
\eqref{mom_bounds}, this implies  the joint limit \eqref{joint_lim}.
\qedhere

\end{proof}

To prove Lemma~\ref{claim}, we use the Laplace transform and
Plancherel argument from the proof of Lemma~\ref{lem_wker}.
Let $\cF_t$ denote the Fourier transform in time with kernel
$e^{i\o t}$, and extend $e^{-t}\wt f_{n,t,0}$ by zero to $t<0$.
Averaging the Laplace transform of the ordered kernel over
permutations before taking the squared spatial norm gives
\[
 \frac{F_n(\o)}{1+\o^2}
 =(2\pi)^n(n!)^2
 \Big\|\cF_t\!\big[e^{-t}\wt f_{n,t,0}\big](\o)
 \Big\|_{\H^{\otimes n}}^2.
\]
That is, with  
$\eta_{\s,k}=\sum_{j=1}^k\xi_{\s(j)}$,
\begin{equation}\label{freq_def}
 F_n(\o)=\int_{\R^n}\prod_{j=1}^n|\xi_j|^{1-2H}
 \bigg|\sum_{\s\in\fS_n}\prod_{k=1}^n
 \frac{1}{(1-i\o)^2+|\eta_{\s,k}|^2}\bigg|^2\,d\pmxi{n}.
\end{equation}
At $\o=0$, the Fourier factor is $1$, so the time transform
reduces to the Laplace transform at parameter $1$. By
\eqref{def_Dn} and \eqref{eq_Y-moment},
$F_n(0)=F_n=(2^H\pi)^n\E[\bfY^n]$.

Lemma~\ref{lem_tail} bounds the integral $\int_\R F_n(\o)/(1+\o^2) d\o$ 
 by $F_n$ times a power of $n$. Lemma~\ref{lem_brown}
then bounds the time-one kernel norm by a moment of $\bfY$.
The power of $n$ drops out after taking $n$th roots, and the
exponential moment bound in Proposition~\ref{prop_exp} gives the
rate in \eqref{ieq_claim}.

\begin{lemma}\label{lem_tail}
Let $-1/4<H<1/2$, and let $\bfY$ and $F_n(\o)$ be defined by \eqref{def_Y} and
\eqref{freq_def}, respectively. The moments of $\bfY$ are given by \eqref{eq_Y-moment}. There exists $C_H>0$ such
that, for every $n\geq1$,
\begin{equation}\label{freq_int}
 \int_\R\frac{F_n(\o)}{1+\o^2}\,d\o
 \leq C_H n^{\frac32-H} \Big(\frac{2\pi}{2^{1-H}}\Big)^n \E [\bfY^n].
\end{equation}
\end{lemma}

\begin{lemma}\label{lem_brown}
Let $-1/4<H<1/2$. There exists $C_H>0$ such that, for all $n\geq1$,
\begin{equation}\label{brown_bd}
 \|\wt f_{n,1,0}\|_{\H^{\otimes n}}^2
 \leq C_H\frac{n^2}{(n!)^{2H+4}}
 \Big(\frac{2^{H-1}}{(1+H)^{2+2H}}\Big)^n\E[\bfY^n].
\end{equation}
\end{lemma}

The proofs of Lemmas~\ref{lem_tail} and~\ref{lem_brown} are given
in Section~\ref{sec_aux}.

\begin{proof}[Proof of Lemma~\ref{claim}]
Fix $0<\ld<\ld_\ast$, where $\ld_\ast = (2/\sfG_H)^{\frac{1+H}{2}}$ is given in \eqref{ieq_exp}.
The term of order $2n$ in the exponential series gives
\[
 \frac{\E[|\bfY|^n]}{(n!)^2} 
 =  \frac{\E[ ( \ld\sqrt{|\bfY|} )^{2n}]}{ \ld^{2n}(n!)^2} 
 \leq \binom{2n}{n}\ld^{-2n}
       \E\!\big[e^{\ld\sqrt{|\bfY|}}\big].
\]
The expectation is finite by Proposition~\ref{prop_exp}, and
$\binom{2n}{n}^{1/n}\to4$ by Stirling's formula.
Taking $n$th roots and then letting
$\ld\uparrow\ld_\ast$, we obtain
\[
 \limsup_{n\to\infty}
 \Big(\frac{\E[|\bfY|^n]}{(n!)^2}\Big)^{1/n}
 \leq\frac{4}{\ld_\ast^2}.
\]
Using this bound with $\E[\bfY^n]\leq\E[|\bfY|^n]$ in \eqref{brown_bd}
and noting that
$(C_H n^2)^{1/n}\to1$, we obtain from \eqref{ieq_exp}
and \eqref{tpHK} that
\begin{align*}
 \limsup_{n\to\infty}
 \big((n!)^{2H+2}\|\wt f_{n,1,0}\|_{\H^{\otimes n}}^2\big)^{1/n}
 &\leq \frac{2^{H-1}}{(1+H)^{2+2H}}\frac{4}{\ld_\ast^2}\\
 &=\frac{\sfG_H^{1+H}}{(1+H)^{2+2H}}
 =\Big(\frac{2\sfK_H}{1+2H}\Big)^{1+2H}.
\end{align*}
Since $1+2H>0$, for each $\dl>0$ this gives \eqref{ieq_claim}
with $C_{H,\dl}=1$ for all sufficiently large $n$. 
Hence, by increasing $C_{H,\dl}$,
the bound  \eqref{ieq_claim} holds true for all $n\geq 1$.
\qedhere

\end{proof}

\subsection{Auxiliary estimates for the upper bound}\label{sec_aux}

\begin{proof}[Proof of Lemma~\ref{lem_tail}]
Let \(F_n\) be given by \eqref{def_Dn}. Concerning \eqref{eq_Y-moment}, it suffices to show 
\begin{equation}\label{freq_int-F}
 \int_\R\frac{F_n(\o)}{1+\o^2}\,d\o
 \leq C_H n^{\frac32-H} F_n.
\end{equation}
Set $q(x)=(1+x^2)^{-1}$ and group the permutations in
\eqref{freq_def} into $\fS_{n,j} \coloneqq \{\s\in\fS_n \colon \s(1)=j\}$,
$1\leq j\leq n$. For a fixed $j$, write $x=\xi_j$ and hold
the other spatial frequencies fixed. For $\s\in\fS_{n,j}$,
the partial sums have the form $\eta_{\s,k}=x+a_{\s,k}$, where
$a_{\s,k}=\sum_{\ell=2}^k\xi_{\s(\ell)}$ is independent of $x$
and $a_{\s,1}=0$. Define
\[
 B_j(x,\o) \coloneqq \sum_{\s\in\fS_{n,j}}\prod_{k=1}^n
 \frac{1}{(1-i\o)^2+(x+a_{\s,k})^2}
 \]
and
\[
A_j(x) \coloneqq B_j(x,0) = \sum_{\s\in\fS_{n,j}}\prod_{k=1}^n
 q(x+a_{\s,k}).
\]
Factoring the denominator gives
$|(1-i\o)^2+u^2|^{-1}=\sqrt{q(u+\o)q(u-\o)}$ for $u\in\R$.
Applying this identity to each factor, the triangle inequality
and Cauchy--Schwarz over $\fS_{n,j}$ yield
\[
 |B_j(x,\o)|^2\leq A_j(x+\o)A_j(x-\o).
\]

To split the spatial weight between the two shifted factors,
set $\rho_n(x)=|x|+n^{-1}$ with
\[
 M_n (x, \o) \coloneqq \rho_n(x+\o) \vee \rho_n(x-\o)
 \AND 
 m_n (x, \o) \coloneqq  \rho_n(x+\o) \wedge \rho_n(x-\o).
\]
The triangle inequality gives
$|x|\leq M_n (x,\o)$
and $|x|\leq |x\pm\o|+|\o|$ for both signs.
Taking the smaller bound and using the fact that $m_n (x,\o) \geq n^{-1}$, we obtain
\[
\begin{aligned}
	|x| &\leq \bigl(|x+\o| \wedge |x-\o|\bigr) + |\o| =m_n (x, \o) -n^{-1} + |\o| \\
 &
 \leq m_n (x, \o) + |\o| \leq(1 + n|\o|) m_n (x, \o).
\end{aligned}
\]
Multiplying $|x|\leq M_n (x, \o)$ by $|x|\leq (1+n|\o|)m_n (x, \o)$ gives
\[
|x|^2\leq (1+n|\o|)M_n (x, \o) m_n (x, \o)
 =(1+n|\o|) \rho_n(x+\o) \rho_n(x-\o).
 \]
Raising both sides to the positive power $\frac12-H$, we obtain
\[
 |x|^{1-2H}
 =\bigl(|x|^2\bigr)^{\frac12-H}
 \leq (1+n|\o|)^{\frac12-H}
       \bigl(\rho_n(x+\o)\rho_n(x-\o)\bigr)^{\frac12-H}.
\]
Combining this with the bound on $|B_j|^2$, applying
Cauchy--Schwarz in $x$, and changing variables $y=x\pm\o$ in
the two integrals, we obtain
\[
 \int_\R |x|^{1-2H}|B_j(x,\o)|^2\,dx
 \leq (1+n|\o|)^{\frac12-H}
       \int_\R \rho_n(y)^{1-2H}A_j(y)^2\,dy.
\]

Since $|(\log q)'(u)|=2|u|/(1+u^2)\leq1$, we have
$q(x+a)\leq e^{|x-y|}q(y+a)$. Each summand of $A_j$ has $n$ such factors,
so $A_j(x)\leq e^{n|x-y|}A_j(y)$.

On each half-line, set $y=x+n^{-1}\sgn(x)$. Then
$|y|=\rho_n(x)$ and $|x-y|=n^{-1}$, so $A_j(x)\leq eA_j(y)$.
The two translations map the half-lines onto $\{y>n^{-1}\}$
and $\{y<-n^{-1}\}$, respectively. Hence
\[
 \int_\R \rho_n(x)^{1-2H}A_j(x)^2\,dx
 \leq e^2\int_{|y|>n^{-1}} |y|^{1-2H}A_j(y)^2\,dy
 \leq e^2\int_\R |y|^{1-2H}A_j(y)^2\,dy.
\]

By Cauchy--Schwarz across the $n$ groups and the preceding
bounds in the variable $x=\xi_j$,
\begin{align*}
 F_n(\o)
 &=\int_{\R^n}\Big|\sum_{j=1}^n B_j(\xi_j,\o)\Big|^2\prod_{\ell=1}^n|\xi_\ell|^{1-2H}\,d\pmxi{n}
 \leq n\sum_{j=1}^n\int_{\R^n}|B_j(\xi_j,\o)|^2\prod_{\ell=1}^n|\xi_\ell|^{1-2H}\,d\pmxi{n}\\
 &\leq e^2 n(1+n|\o|)^{\frac12-H}
       \sum_{j=1}^n\int_{\R^n}A_j(\xi_j)^2\prod_{\ell=1}^n|\xi_\ell|^{1-2H}\,d\pmxi{n}
       \leq e^2 n(1+n|\o|)^{\frac12-H}F_n.
\end{align*}
The last inequality uses the nonnegativity of the group sums
at $\o=0$, whose sum is the permutation sum defining $F_n$.
Since $1+n|\o|\leq n(1+|\o|)$ and $\frac12-H<1$, 
we have
\[
 \int_\R\frac{F_n(\o)}{1+\o^2}\,d\o
 \leq e^2 n^{\frac32-H}F_n
       \int_\R\frac{(1+|\o|)^{\frac12-H}}{1+\o^2}\,d\o
 \leq C_H n^{\frac32-H}F_n.
\]
The proof of \eqref{freq_int-F} and hence Lemma \ref{lem_tail} is complete. 
\end{proof}

\begin{proof}[Proof of Lemma~\ref{lem_brown}]
Set $b=2H+2$. Integrating the scaling identity \eqref{coeff_scale}
against $e^{-2t}\,dt$ gives the first equality below. The second
follows from Plancherel \eqref{plan} in time and the Fourier
representation of $F_n(\o)/(1+\o^2)$ preceding \eqref{freq_def}.
\begin{equation}\label{gamma_id}
\begin{aligned}
 &\frac{\Gam(bn+1)}{2^{bn+1}}\|\wt f_{n,1,0}\|_{\H^{\otimes n}}^2
 =\int_0^\infty e^{-2t}\|\wt f_{n,t,0}\|_{\H^{\otimes n}}^2\,dt\\
 & \qquad =\frac{1}{(2\pi)^{n+1}(n!)^2}
   \int_\R\frac{F_n(\o)}{1+\o^2}\,d\o
   \;\leq C_H\frac{n^{\frac32-H}}{(n!)^2}
   2^{(H-1)n}\E[\bfY^n],
\end{aligned}
\end{equation}
where the inequality follows from \eqref{freq_int} and
$(2^H\pi)^n/(2\pi)^{n+1}=2^{(H-1)n}/(2\pi)$.
For analogous comparisons when the spatial
covariance has a nonnegative convolution square root, see
\cite[Lemmas~3.5 and~3.7]{BCC22}.

By Stirling's formula \eqref{stir},
\[
 \frac{1}{\Gam(bn+1)}
 \leq C_H\frac{n^{(b-1)/2}}{b^{bn}(n!)^b}.
\]
Multiplying \eqref{gamma_id} by $2^{bn+1}/\Gam(bn+1)$ and
applying this estimate gives
\[
 \|\wt f_{n,1,0}\|_{\H^{\otimes n}}^2
 \leq C_H\frac{n^{\frac32-H+(b-1)/2}}{(n!)^{b+2}}
 \left[2^{H-1}\left(\frac{2}{b}\right)^b\right]^n\E[\bfY^n].
\]
Since $b=2+2H$, the exponent of $n$ equals $2$ and
$2^{H-1}(2/b)^b=2^{H-1}/(1+H)^{2+2H}$. This proves \eqref{brown_bd}.
\end{proof}

\section*{Development of the work and use of AI}
\small This work originated during Xiong Wang's visit to Boston University on December 6--9, 2024. At the SPDE workshop held at the Brin Mathematics Research Center in June 2025, the authors established the existence and uniqueness of the random-field solution for \(- 1/4 < H < 1/2\), as well as the moment upper bound stated in Proposition~\ref{prop_ham_up}.  At that stage, the authors had obtained two-sided
moment bounds and had not yet to identify the exponential growth constant
$\sfK_H$ in \eqref{joint_lim}.

In  September--October of 2026, the authors used OpenAI Codex to explore the
exact upper bound under the authors' prompts without success on the whole region $-1/4<H<1/2$. 
However, Codex proposed applying the compactification of
Brownian occupation measures developed by Mukherjee and Varadhan
\cite{MV16} and produced an argument for the exponential moment
estimate in Proposition~\ref{prop_exp}. This estimate supplies
the sharp constant in the upper bound.
Also, when the authors proposed using Schr\"odinger operators to prove the
exact lower bound in Proposition~\ref{prop_mom}, Codex identified the
relevant spectral results in Teschl's book~\cite{Tes14}.
During revision,  the authors spent substantial time working
through these arguments and used GPT-6.1 Astra to examine individual
steps and refine the exposition. The authors take responsibility
for the mathematical arguments and the manuscript in its present
form.


\begin{thebibliography}{99}



\bibitem[Bal12]{Bal12}
R.M.~Balan,
{\it The stochastic wave equation with multiplicative fractional noise: a Malliavin calculus approach,}
Potential Anal. \textbf{36} (2012), no. 1, 1--34.

\bibitem[BCC22]{BCC22}
R.M.~Balan, L.~Chen, and X.~Chen,
{\it Exact asymptotics of the stochastic wave equation with
time-independent noise,}
Ann. Inst. Henri Poincar\'e Probab. Stat. \textbf{58} (2022), no.~3,
1590--1620.

\bibitem[BCM22]{BCM22}
R.M.~Balan, L.~Chen, and Y.~Ma,
{\it Parabolic Anderson model with rough noise in space and rough
initial conditions,}
Electron. Commun. Probab. \textbf{27} (2022), Paper No.~65, 12 pp.


\bibitem[BJQS15]{BJQS15}
R.M.~Balan, M.~Jolis, and L.~Quer-Sardanyons,
{\it SPDEs with affine multiplicative fractional noise in space with index $\frac{1}{4} < H < \frac{1}{2}$,}
Electron. J. Probab. \textbf{20} (2015), no. 54, 36 pp.

\bibitem[BJQS17]{BJQS17}
R.M.~Balan, M.~Jolis, and L.~Quer-Sardanyons,
{\it Intermittency for the hyperbolic Anderson model with rough noise in
space,}
Stochastic Process. Appl. \textbf{127} (2017), no.~7, 2316--2338.

\bibitem[BNQZ22]{BNQZ22}
R.M.~Balan, D.~Nualart, L.~Quer-Sardanyons, and G.~Zheng,
{\it The hyperbolic Anderson model: moment estimates of the Malliavin
derivatives and applications,}
Stoch. Partial Differ. Equ. Anal. Comput. \textbf{10} (2022), no.~3,
757--827.

\bibitem[BS17]{BS17}
R.M.~Balan and J.~Song,
{\it Hyperbolic Anderson model with space-time homogeneous Gaussian
noise,}
ALEA Lat. Am. J. Probab. Math. Stat. \textbf{14} (2017), no.~2,
799--849.

\bibitem[BS19]{BS19}
R.M.~Balan and J.~Song,
{\it Second order Lyapunov exponents for parabolic and hyperbolic
Anderson models,}
Bernoulli \textbf{25} (2019), no.~4A, 3069--3089.

\bibitem[BY24]{BY24}
R.M.~Balan and W.~Yuan,
{\it Hyperbolic Anderson model with time-independent rough noise:
Gaussian fluctuations,}
Electron. J. Probab. \textbf{29} (2024), Paper No.~177, 52 pp.


\bibitem[CM94]{CM94}
R.A.~Carmona and S.A.~Molchanov,
{\it Parabolic Anderson problem and intermittency,}
Mem. Amer. Math. Soc. \textbf{108} (1994), no.~518, viii+125 pp.

\bibitem[Che26]{Che26}
L.~Chen,
{\it Stratonovich waves driven by space-time Gaussian noise,}
arXiv:2611.23014 (2026).

\bibitem[CD15a]{CD15-1}
L.~Chen and R.C.~Dalang,
{\it Moments and growth indices for the nonlinear stochastic heat equation with rough initial conditions,}
Ann. Probab. \textbf{43} (2015), no. 6, 3006--3051.

\bibitem[CD15b]{CD15-2}
L.~Chen and R.C.~Dalang,
{\it Moment bounds and asymptotics for the stochastic wave equation,}
Stochastic Process. Appl. \textbf{125} (2015), no. 4, 1605--1628.


\bibitem[CLX26]{CLX26}
L.~Chen, C.Y.~Lee and P.~Xia,
{\it Strong local nondeterminism for stochastic time-fractional slow and fast diffusion equations,}
arXiv preprint \href{https://arxiv.org/abs/2602.05317}{arXiv:2602.05317} (2026).

\bibitem[Che07]{Che07}
X.~Chen,
{\it Large deviations and laws of the iterated logarithm for the
local times of additive stable processes,}
Ann. Probab. \textbf{35} (2007), no.~2, 602--648.


\bibitem[Che16]{Che16}
X.~Chen,
{\it Spatial asymptotics for the parabolic Anderson models with generalized time-space Gaussian noise,}
Ann. Probab. \textbf{44} (2016), no. 2, 1535--1598.

\bibitem[Che17]{Che17}
X.~Chen,
{\it Moment asymptotics for parabolic {A}nderson equation with fractional time-space noise: in {S}korokhod regime,}
Ann. Inst. Henri Poincar\'e Probab. Stat. \textbf{53} (2017), no.~2,
819--841.

\bibitem[Che19]{Che19}
X.~Chen,
{\it Parabolic Anderson model with rough or critical Gaussian noise,}
Ann. Inst. Henri Poincar\'e Probab. Stat. \textbf{55} (2019), no.~2,
941--976.

\bibitem[Che23]{Chen23SILT}
X.~Chen,
{\it Exponential asymptotics for Brownian self-intersection
local times under Dalang's condition,}
Electron. J. Probab. \textbf{28} (2023), paper no.~119, 1--17.


\bibitem[CDOT21]{CDOT21}
X.~Chen, A.~Deya, C.~Ouyang, and S.~Tindel,
{\it Moment estimates for some renormalized parabolic Anderson models,}
Ann. Probab. \textbf{49} (2021), no. 5, 2599--2636.

\bibitem[CDST23]{CDST23}
X.~Chen, A.~Deya, J.~Song, and S.~Tindel,
{\it Hyperbolic Anderson model 2: Strichartz estimates and Stratonovich setting,}
Int. Math. Res. Not. IMRN (2023), no.~21,
18575--18628.


\bibitem[CDST25]{CDST25}
X.~Chen, A.~Deya, J.~Song, and S.~Tindel,
{\it Solving the hyperbolic Anderson model 1: Skorohod setting,}
Ann. Inst. Henri Poincar\'e Probab. Stat. \textbf{61} (2025), no.~3,
1794--1814.

\bibitem[CH25]{CH25}
X.~Chen and Y.~Hu,
{\it Hyperbolic Anderson equations with general time-independent
Gaussian noise: Stratonovich regime,}
Ann. Probab. \textbf{53} (2025), no.~6, 2144--2195.


\bibitem[CHSX15]{CHSX15}
X.~Chen, Y.~Hu, J.~Song, and F.~Xing,
{\it Exponential asymptotics for time-space Hamiltonians,}
Ann. Inst. Henri Poincar\'e Probab. Stat. \textbf{51} (2015), no.~4,
1529--1561.


\bibitem[CH23]{CH23}
Z.~Chen and Y.~Hu,
{\it Solvability of parabolic Anderson equation with fractional Gaussian noise,}
Commun. Math. Stat. \textbf{11} (2023), no. 3, 563--582.

\bibitem[Dal99]{Dal99}
R.C.~Dalang,
{\it Extending the martingale measure stochastic integral with
applications to spatially homogeneous S.P.D.E.'s,}
Electron. J. Probab. \textbf{4} (1999), no.~6, 29 pp.

\bibitem[DM09]{DM09}
R.C.~Dalang and C.~Mueller,
{\it Intermittency properties in a hyperbolic Anderson problem,}
Ann. Inst. Henri Poincar\'e Probab. Stat. \textbf{45} (2009), no.~4,
1150--1164.

\bibitem[DZ10]{DZ10}
A.~Dembo and O.~Zeitouni,
{\it Large Deviations Techniques and Applications,}
2nd ed., corrected printing, Stochastic Modelling and Applied
 Probability,
vol.~38, Springer, Berlin, 2010.

\bibitem[DPV12]{DPV12}
E.~Di Nezza, G.~Palatucci, and E.~Valdinoci,
{\it Hitchhiker's guide to the fractional Sobolev spaces,}
Bull. Sci. Math. \textbf{136} (2012), no.~5, 521--573.

\bibitem[Ger12]{Ger12}
S.~Gerhold,
{\it Asymptotics for a variant of the Mittag--Leffler function,}
Integral Transforms Spec. Funct. \textbf{23} (2012), no.~6, 397--403.

\bibitem[GSW26]{GSW26}
Y.~Guo, J.~Song, and G.~Wu,
{\it Stochastic fractional diffusion equations with spatial Gaussian noise in the Stratonovich regime,}
arXiv preprint \href{https://arxiv.org/abs/2609.04826}{arXiv:2609.04826} (2026).

\bibitem[Hu02]{Hu02}
Y.~Hu,
{\it Chaos expansion of heat equations with white noise potentials,}
Potential Anal. \textbf{16} (2002), no. 1, 45--66.

\bibitem[Hu17]{Hu17}
Y.~Hu,
{\it Analysis on Gaussian spaces,} World Scientific Publishing Co. Pte. Ltd., Hackensack, NJ, 2017. 

\bibitem[HHLNT17]{HHLNT17}
Y.~Hu, J.~Huang, K.~L\^e, D.~Nualart, and S.~Tindel,
{\it Stochastic heat equation with rough dependence in space,}
Ann. Probab. \textbf{45} (2017), no. 6B, 4561--4616.

\bibitem[HHNT15]{HHNT15}
Y.~Hu, J.~Huang, D.~Nualart, and S.~Tindel,
{\it Stochastic heat equations with general multiplicative Gaussian
noises: H\"older continuity and intermittency,}
Electron. J. Probab. \textbf{20} (2015), Paper No.~55.

 
\bibitem[Kho14]{Kho14}
D.~Khoshnevisan,
{\it Analysis of stochastic partial differential equations,}
CBMS Regional Conference Series in Mathematics, vol.~119,
American Mathematical Society, Providence, RI, 2014.

\bibitem[Kon16]{Kon16}
W.~K\"onig,
{\it The parabolic {A}nderson model: Random walk in random potential,}
Birkh\"auser/Springer, Cham, 2016.

\bibitem[LHW24]{LHW24}
S.~Liu, Y.~Hu, and X.~Wang,
{\it In search of necessary and sufficient conditions to solve the
parabolic Anderson model with fractional Gaussian noises,}
Electron. J. Probab. \textbf{29} (2024), Paper No.~140, 48 pp.

\bibitem[LHW26]{LHW26}
S.~Liu, Y.~Hu, and X.~Wang,
{\it Stochastic wave equation with additive fractional noise: solvability and global H\"older continuity,}
Stochastic Process. Appl. {\bf 198} (2026), Paper No. 104964, 24 pp.

\bibitem[LSSW16]{LSSW16}
A.~Lodhia, S.~Sheffield, X.~Sun, and S.S.~Watson,
{\it Fractional Gaussian fields: A survey,}
Probab. Surv. \textbf{13} (2016), 1--56.

\bibitem[MV16]{MV16}
C.~Mukherjee and S.R.S.~Varadhan,
{\it Brownian occupation measures, compactness and large deviations,}
Ann. Probab. \textbf{44} (2016), no.~6, 3934--3964.

\bibitem[NP12]{NP12}
I.~Nourdin and G.~Peccati,
{\it Normal approximations with Malliavin calculus: From Stein's method
to universality,}
Cambridge Tracts in Mathematics, vol.~192,
Cambridge University Press, Cambridge, 2012.

\bibitem[Nua06]{Nua06}
D.~Nualart,
{\it The Malliavin calculus and related topics,}
second edition, Springer-Verlag, Berlin, 2006.

\bibitem[Sim82]{Sim82}
B.~Simon,
{\it Schr\"odinger semigroups,}
Bull. Amer. Math. Soc. (N.S.) \textbf{7} (1982), no.~3, 447--526.

\bibitem[SSX20]{SSX20}
J.~Song, X.~Song, and F.~Xu,
{\it Fractional stochastic wave equation driven by a Gaussian noise
rough in space,}
Bernoulli \textbf{26} (2020), no.~4, 2699--2726.

\bibitem[Tes14]{Tes14}
G.~Teschl,
{\it Mathematical methods in quantum mechanics: With applications to
Schr\"odinger operators,}
second edition, Graduate Studies in Mathematics, vol.~157,
American Mathematical Society, Providence, RI, 2014.

\bibitem[Wal86]{Wal86}
J.B.~Walsh,
{\it An introduction to stochastic partial differential equations,}
in \emph{\'Ecole d'\'Et\'e de Probabilit\'es de Saint-Flour XIV--1984},
Lecture Notes in Math., vol.~1180, Springer, Berlin, 1986, pp.~265--439.

\end{thebibliography}
\end{document}